\documentclass{article}
\usepackage{verbatim,amsmath,amssymb,mathrsfs,amsthm}
\usepackage{fancyhdr}
\usepackage[dvipsnames]{xcolor}
\usepackage{graphicx,wrapfig}
\usepackage{subfigure}
\usepackage{caption}
\usepackage{amsmath}
\usepackage{makecell}
\usepackage{multirow}
\usepackage{array}
\usepackage{setspace}
\usepackage{float}
\usepackage{color,mathtools}
\usepackage{epstopdf}
\usepackage{lipsum}
\usepackage{appendix}
\usepackage{algorithm}
\usepackage{algorithmic}
\usepackage{blindtext}
\usepackage[textwidth=15.8cm]{geometry}

\newcommand{\bs}[1]{\boldsymbol{#1}}
\newcommand{\bmu}{\bs{\mu}}

\newcommand{\bc}{\bs{c}}

\newcommand{\calN}{\mathcal{N}}
\newcommand{\calD}{\mathcal{D}}
\newcommand{\N}{{\cal N}}
\DeclareMathOperator*{\argmax}{arg\,max}
\DeclareMathOperator*{\argmin}{arg\,min}

\newcommand{\bx}{\bs{x}}

\usepackage{amsthm,amsmath,amssymb}
\newtheorem{theorem}{Theorem}

\newtheorem{remark}{Remark}


\begin{document}

\graphicspath{{Figs/}}

\title{Fast $L^2$ optimal mass transport via reduced basis methods for the Monge-Amp$\grave{\rm e}$re equation}
\author{Shijin Hou\thanks{School of Mathematical Sciences, University of Science and Technology of China, Hefei, Anhui
230026, People's Republic of China. Email: {\tt{houshiji@mail.ustc.edu.cn}}.}
\and Yanlai Chen\thanks{Department of Mathematics, University of Massachusetts Dartmouth, 285 Old Westport Road, North Dartmouth, MA 02747, USA. Email: {\tt{yanlai.chen@umassd.edu}}. This author was partially supported by National Science Foundation grant DMS-1719698, by a grant from the College of Arts \& Sciences at the University of Massachusetts Dartmouth, and by the UMass Dartmouth Marine and UnderSea Technology (MUST) Research Program made possible via an Office of Naval Research grant N00014-20-1-2849.}
\and Yinhua Xia\thanks{School of Mathematical Sciences, University of Science and Technology of China, Hefei, Anhui
	230026, People's Republic of China. Email: {\tt{yhxia@ustc.edu.cn}}. This author was partially supported by National Natural Science Foundation of China  grant No. 11871449}}
\date{}

\maketitle

\begin{abstract}
Repeatedly solving the parameterized optimal mass transport (pOMT) problem is a frequent task in applications such as image registration and adaptive grid generation. It is thus critical to develop a highly efficient reduced solver that is equally accurate as the full order model. In this paper, we propose such a machine learning-like method for pOMT by adapting a new reduced basis (RB) technique specifically designed for nonlinear equations, the reduced residual reduced over-collocation (R2-ROC) approach, to the parameterized Monge-Amp$\grave{\rm e}$re equation. It builds on top of a narrow-stencil finite different method (FDM), a so-called truth solver, which we propose in this paper for the Monge-Amp$\grave{\rm e}$re equation with a transport boundary.  Together with the R2-ROC approach, it allows us to handle the strong and unique nonlinearity pertaining to the Monge-Amp$\grave{\rm e}$re equation achieving online efficiency without resorting to any direct approximation of the nonlinearity. Several challenging numerical tests demonstrate the accuracy and high efficiency of our method for solving the Monge-Amp$\grave{\rm e}$re equation with various parametric boundary conditions.
\end{abstract}

\section{Introduction}
The Optimal Mass Transport (OMT) problem has received significant attention in recent years thanks to its wide applicability in areas such as image retrieval \cite{li2013novel,rubner2000earth}, shape and image registration \cite{haker2004optimal,lai2014multi}, super-resolution reconstruction \cite{kolouri2015transport}, cancer detection \cite{basu2014detecting,ozolek2014accurate}, machine
learning \cite{Goodfellow2014,Arjovsky2017,el2012bayesian,frogner2015learning,kolouri2016sliced,montavon2015wasserstein} and adaptive grid generation \cite{sulman2011optimal}, just to name a few.
Among these applications, an interesting scenario emerges when the OMT problem needs to be solved repeatedly and often in a real-time manner. For example, in image processing, solving a OMT problem provides the optimal transformation between one pair of images out of potentially many that are closely related. Another example is that one OMT problem needs to be resolved for determining the grid movement in adaptive grid generation \cite{sulman2011optimal} for every round of a posterior error estimation. Given appropriate parameterizations, these problems can be modeled by a pOMT problem, the focus of the current paper.

Initially proposed by Gaspard Monge \cite{monge1781memoire} in the 18$^{\rm th}$ century, OMT seeks the optimal mass-preserving transform between two distributions of mass for a given cost of transportation.
Given two bounded and open domains $X,Y\in \mathbb{R}^{d}$, let $\nu_X(\pmb{x}, \bmu)$ be a probability measure on $X$, parameterized by a $p$-dimensional parameter $\bmu \in \mathcal{D}\subset \mathbb{R}^{p}$. $\nu_Y(\pmb{y}, \bmu) = T_{\#} \nu_X$ is its push-forward on $Y$ with a measurable map $T: X \mapsto Y$ satisfying
\[
\forall h \in {\mathcal C}(Y) \quad \int_Y h(\pmb{y}) d \nu_Y(\pmb{y}, \bmu) = \int_X h(T(\pmb{x})) d \nu_X(\pmb{x}, \bmu) .
\]
OMT seeks a minimizer of the cost functional
$I(T) = \int_{X}c(\pmb{x},T(\pmb{x}))d\nu_X(\pmb{x}, \bmu)$,
where $c(\pmb{x},\pmb{y})$ denotes the cost of transporting a unit of mass from the point $\pmb{x}\in X$ to the point $\pmb{y}\in Y$. If the measures are absolutely continuous with (parametric) positive densities $f_{X}(\pmb{x}, \bmu), f_{y}(\pmb{y},\bmu)$, that is,
\[
d\nu_X(\pmb{x}, \bmu) = f_X(\pmb{x}, \bmu)d\pmb{x},\ d\nu_Y(\pmb{y}, \bmu) = f_Y(\pmb{y}, \bmu)d\pmb{y},
\]
by simple calculation, the mass balance condition could be rewritten as
\begin{equation}\label{eq_mass_conserved}
{\rm det}\left(DT(\pmb{x})\right)f_{Y}(T(\pmb{x}), \bmu) = f_{X}(\pmb{x}, \bmu),
\end{equation}
where ${\rm det}\left(DT(\pmb{x})\right)$ denotes the determinant of the Jacobian of $T(\pmb{x})$. Although there are other formulations of this problem such as Kantorovich formulation \cite{kantorovich2006problem}, we focus on the Monge formulation in this paper and aim to develop a fast solver for it. In the special case of the quadratic cost function
$c(\pmb{x},\pmb{y}) = \frac{1}{2}\vert\pmb{x}-\pmb{y}\vert^{2}$,
the minimizing map $T(\pmb{x})$ can be expressed as the gradient of a convex function \cite{evans1997partial,rockafellar1966characterization} justifying a substitution of $T$ by $\nabla u$ in (\ref{eq_mass_conserved}). This results in the following parametric Monge-Amp$\grave{\rm e}$re equation which is augmented by the convexity constraint on $u$ for uniqueness and stability  \cite{aleksandrov1960certain,gutierrez2001monge}. We also enforce the so-called transport or second boundary condition.
\begin{subequations}
\label{eq_MA}
\begin{align}
\label{eq_MA_with_transform_BC}
&{\rm det}\left(D^{2}u(\pmb{x}, \bmu)\right)  =\frac{f_{X}(\pmb{x}, \bmu)}{f_{Y}(\nabla u(\pmb{x},\bmu), \bmu)},\  \pmb{x}\in X,\\
\label{BC_transport}
&\nabla u(X, \bmu)  = Y,\\
\label{MA_convexity}
&u(\pmb{x}, \bmu)\ \rm{is\ convex\ in}\ X.
\end{align}
\end{subequations}
Here, $D^{2}u(x)$ denotes the Hessian of the function $u$. The challenge in solving \eqref{eq_MA} resides in the strong nonlinearity, the convexity constraint, the difficulty of approximating the transport boundary condition, and the low regularity of its solution. The literature of numerical methods is therefore rather scarce. Benamou and Brenier \cite{benamou2000computational} presented a fluid flow approach which was further developed by Haber et. al. \cite{haber2010efficient}. The method
is computationally expensive due to the introduction of an additional dimension. More recently, Froese \cite{froese2012numerical} proposed an approach for solving \eqref{eq_MA} by iteratively solving
a sequence of Neumann boundary value problems of the Monge-Amp$\grave{\rm e}$re equation and proved that the numerical solution converges to the viscosity solution. Benamou et. al. \cite{benamou2012viscosity} developed an approach by reformulating the transport boundary condition by an implicit Hamilton-Jacobi equation and gave the proof of convergence.

The situation is exacerbated by the need to solve the problem for a large number of parameter values. To the best of our knowledge, there is no existing work based on systematic model order reduction for the parameterized Monge-Amp$\grave{\rm e}$re equation. In this article, we aim to provide {\em first such work}. Our first contribution is the proposal of a truth solver by extending the narrow-stencil finite difference scheme of \cite{feng2019narrow}, originally designed for the Hamilton-Jacobi-Bellman equations, to our context.
An improvement of the standard finite difference scheme, this narrow-stencil scheme is amenable to the RB framework while being more robust in handling singular solutions thanks to the introduction of the artificial viscosity and numerical moment. In addition, we adopt the framework of \cite{froese2012numerical} in dealing with the transport boundary. The R2-ROC method proposed in \cite{chen2021eim} serves as our reduced order modeling approach. The R2-ROC is a class of Reduced Basis Method (RBM) \cite{Quarteroni2015, HesthavenRozzaStammBook, Haasdonk2017Review} specifically designed for nonlinear and nonaffine problems. Like RBM but based on an underlying scheme of a nodal form, it features an offline/online decomposition strategy, a posteriori error estimator/indicator, and a classical greedy algorithm. The main task of the offline phase is to construct a problem-dependent, low-dimensional surrogate space and set the stage for the online computations.
After the (time consuming) offline stage, the full speed of the method will then be appreciated online when the reduced solver is performed on demand and usually with a cost only dependent on the (much lower) RB space dimension. Due to the strong nonlinearity of the Monge-Amp$\grave{\rm e}$re equation, the classical RBM will suffer on its online complexity resulting from its dependence on the number of EIM/DEIM decomposition \cite{barrault2004empirical,chaturantabut2010nonlinear,grepl2007efficient} terms. The R2-ROC method eliminates this dependence by augmenting and extending the EIM approach as a direct PDE solver, judiciously determining a set of over-collocation points, and taking advantage of the simplicity of evaluating the hyper-reduced well-chosen residuals. It achieves offline/online computation efficiency and, more interestingly, the independence on the number of EIM/GEIM expansion terms.
Our second contribution of this paper is to adapt the R2-ROC, designed for the classical Dirichlet or Neumann boundary conditions, for the much more intricate transport boundary condition imposed by the Monge-Amp$\grave{\rm e}$re equation.

The organization of this paper is as follows. In Section \ref{Section:fdm}, we review some theoretical results for the Monge-Amp$\grave{\rm e}$re equation before describing an iterative algorithm for implementing the transport boundary condition and a narrow-stencil finite difference scheme for approximating the Monge-Amp$\grave{\rm e}$re equation. The combination of them provides an efficient full order model for the transport boundary problem of the Monge-Amp$\grave{\rm e}$re equation. In Section \ref{sec:R2-ROC-Alg}, we introduce the R2-ROC method and our adaptation to the transport boundary case toward our reduced order model. In Section \ref{Section:numerical}, we present numerical results for the Monge-Amp$\grave{\rm e}$re equation and the parameterized Monge-Amp$\grave{\rm e}$re equation to demonstrate the efficiency and accuracy of our methods. Concluding remarks are made in Section \ref{Section:conclusion}.

\section{A narrow-stencil finite difference method for $L^{2}$ optimal mass transport problem}
\label{Section:fdm}

This section is devoted to a detailed description of our truth solver, an extension of the narrow-stencil finite difference scheme of \cite{feng2019narrow} to the Monge-Amp$\grave{\rm e}$re equation adopting the framework of \cite{froese2012numerical} in dealing with the transport boundary condition. In order to properly inform the numerical scheme, it is important to know when a classical $C^{2}$ solution exists. Toward that end, we first review some regularity results for the Monge-Amp$\grave{\rm e}$re equation.

\subsection{Regularity}
\label{regularity}
The classical $C^{2}$ solution of the Monge-Amp$\grave{\rm e}$re type equation exist under certain regularity conditions on the data and computational domains. We present below the regularity results for the Dirichlet boundary value problem
\begin{subequations}
\label{eq_MA_D}
\begin{align}
&{\rm det}(D^{2}u(\pmb{x}))=f(\pmb{x}),\  \pmb{x}\in X,\label{eq_MA_with_Dirichlet_BC}\\
&u(\pmb{x}) = g(\pmb{x}),\ \pmb{x}\in \partial X.\label{BC_Dirichlet},
\end{align}
\end{subequations}
and then the transport boundary value problem \eqref{eq_MA}.
\begin{theorem}\label{regularity_1}
\cite{Urbas1986The,caffarelli1984dirichlet}
Suppose that X is strictly convex with boundary $\partial X \in C^{2,\alpha}$. Suppose also that the function $f\in C^{\alpha}(X)$ is strictly positive and the boundary values $g\in  C^{2,\alpha}(\partial X)$. Then the Dirichlet boundary value problem of Monge-Amp$\grave{e}$re equation \eqref{eq_MA_D} has a unique $C^{2,\alpha}$ solution.
\end{theorem}

\begin{theorem}\label{regularity_2}
(Interior Regularity \cite{caffarelli1991some,Caffarelli1992The}) Suppose that X, Y are bounded, connected, open sets and Y is convex. Suppose also that the density functions
\begin{equation*}
f_{X}: X \to (0,\infty),\ f_{Y}: Y\to (0,\infty)
\end{equation*}
are bounded away from 0 and $\infty$. Then the solution of the Monge-Amp$\grave{e}$re equation \eqref{eq_MA} belongs to $C^{1,\alpha}(X)$ for some $0<\alpha<1$.
If , in addition, the density function $f,\ g\in C^{\beta}$ for some $0<\beta<1$ then the solution of the Monge-Amp$\grave{e}$re equation belongs to $C^{2,\alpha}$ for every $0<\alpha<\beta$.
\end{theorem}

\subsection{Transport boundary condition}\label{Section:Transport boundary condition}

In this section, we describe an efficient algorithm proposed in \cite{froese2012numerical} for dealing with the transport condition \eqref{BC_transport}. Indeed, when $X,\ Y$ are both convex, the transport condition can be enforced by requiring that the boundary points of $X$ are mapped to the boundary points of $Y$ \cite{pogorelev1971dirichlet,trudinger2009second,urbas1997second}. That is,
\[
\nabla u(\pmb{\mu}):\partial X\mapsto \partial Y.
\]
Assuming that the boundary of target region $\partial Y$ can be represented as the zero$^{\rm th}$ level set of a function $\Phi$, we have that the transport map must satisfy
\begin{equation}\label{BC_transport_exact}
\Phi(\nabla u(\pmb{x};\pmb{\mu})) = 0,\ \pmb{x}\in \partial X.
\end{equation}
The appearance of the gradient and the simplicity of implementing a Neumann boundary condition motivated the authors of \cite{froese2012numerical} to replace the condition (\ref{BC_transport_exact}) by the Neumann boundary condition
\[
\nabla u(\pmb{x};\pmb{\mu})\cdot \pmb{n}(\pmb{x})=\phi(\pmb{x};\pmb{\mu}),
\]
where $\phi(\pmb{x};\pmb{\mu})$ is a function to be determined and $\pmb{n}$ denotes the unit outward normal.
\cite{froese2012numerical}  further proposes the following iterative approach for
approximating the function $\phi(\pmb{x};\pmb{\mu})$. Given the $k^{\rm th}$ iterate $u^{k}(\pmb{x};\pmb{\mu})$ of the approximate solution to the Monge-Amp$\grave{\rm e}$re equation with transport boundary condition, we proceed to the next iterate as follows.
We first compute $\phi^{k}(\pmb{x};\pmb{\mu})$ for $\pmb{x}\in \partial X$ via
\begin{equation}
\label{eq_Neumann_BC_of_subproblem}
\phi^{k}(\pmb{x};\pmb{\mu}) = \mathbb{P}_{\partial Y}(\nabla u^{k}(\pmb{x};\pmb{\mu}))\cdot \pmb{n}(\pmb{x}),
\end{equation}
where $\mathbb{P}_{\partial Y}(\pmb{v})$ denotes the shortest-distance projection of $\pmb{v}$ onto the set $\partial Y$:
$\mathbb{P}_{\partial Y}(\pmb{v})=\mathop{\arg\min}_{\pmb{w}\in \partial Y}||\pmb{w}-\pmb{v}||_{L^2}^{2}.$
After that, we find a convex function $u^{k+1}(\pmb{x};\pmb{\mu}): X\mapsto \mathbb{R}$ and a constant $\sigma^{k+1}(\pmb{\mu})\in R^{+}$ such that
\begin{subequations}
\label{eq_subproblem}
\begin{align}
&{\rm det}(D^{2}u^{k+1}(\pmb{x};\pmb{\mu})) = \sigma^{k+1}(\pmb{\mu})\frac{f_{X}(\pmb{x};\pmb{\mu})}{\widehat{f_{Y}}(\nabla u^{k+1}(\pmb{x};\pmb{\mu});\pmb{\mu})}
\label{eq_subproblem_a}
\\%
&\int_{X}u^{k+1}(\pmb{x};\pmb{\mu})dx =0,\label{eq_subproblem_b} \\
&\nabla u^{k+1}(\pmb{x};\pmb{\mu})\cdot \pmb{n}(\pmb{x}) = \phi^{k}(\pmb{x};\pmb{\mu}),\ \pmb{x}\in\partial X.
\label{eq_subproblem_c}.
\end{align}
\end{subequations}
Here, $\widehat{f_{Y}}(\pmb{y};\pmb{\mu})$ is the extended target density function of $f_{Y}(\pmb{y};\pmb{\mu})$ defined as
\[
\widehat{f_{Y}}(\pmb{y};\pmb{\mu})=
\left\{
\begin{aligned}
&f_{Y}(\pmb{y};\pmb{\mu}),\ \pmb{y}\in Y,\\
&f_{Y}(\pmb{y}_{0};\pmb{\mu}),\ \pmb{y}\notin Y,
\end{aligned}
\right.
\]
where $\pmb{y}_0$ is a point in the interior of the target set $Y$. This extension assigns positive values outside of $Y$ to accommodate the fact that it is a density and, during iterations, $\nabla u^{k+1}$ in \eqref{eq_subproblem_a} may map (part of) $X$ out of $Y$.  It's worth noting that more complex extensions with higher regularity may be required for the convergence of the Monge-Amp$\grave{\rm e}$re solver, see \cite{froese2012numerical}. We adopt this simple position extension of $f_{Y}(\pmb{y};\pmb{\mu})$ which is enough for our solver introduced in Section \ref{Section:The narrow-stencil finite difference scheme}. This iteration proceeds until the difference between $u^{k}$ and $u^{k+1}$ is sufficiently small. To start it, we simply let $\phi^{-1}(\pmb{x};\pmb{\mu}) = B\pmb{x}\cdot \pmb{n}$, where $B>0$ is large enough to ensure the set $\{B\pmb{x}\ |\ \pmb{x}\in X\}$ contains $Y$.

\begin{remark}
As is well known, the solution of the Neumann boundary value problems may not exist. Even if it exists, the solution is unique only up to a constant. For these reasons, the variable $\sigma^{k+1}(\pmb{\mu})$ in \eqref{eq_subproblem_a} and the mean-zero condition \eqref{eq_subproblem_b} are introduced. The projection $\mathbb{P}_{\partial Y}(\pmb{v})$ is introduced to mitigate the misalignment of $\nabla u^{k}(\partial X;\pmb{\mu})$ and $\partial Y$, which contributes to obtain the correct choice of $\phi(\pmb{x};\pmb{\mu})$.
\end{remark}
\begin{remark}
An $L^{2}$ optimal mapping does not lead to twisting or rotation. When we consider the simple case of mapping a rectangle to another rectangle, each side of $X$ will be mapped to the corresponding sides of $Y$. Since the directional derivative of $u$ at each $\pmb{x}\in \partial X$ is determined, we obtain an exact Neumann boundary condition. In this case, the transport boundary problem becomes a Neumann boundary problem.
\end{remark}

\subsection{Discretization for the Monge-Amp$\grave{\rm \textbf{e}}$re equation }\label{Section:The narrow-stencil finite difference scheme}
With the iterative framework for the transport problem, we propose to incorporate a finite difference solver \cite{lewis2013finite, feng2019narrow} for the Neumann boundary value problems \eqref{eq_subproblem} which would conclude the description of our full order model.
This finite difference solver adopts artificial viscosity and moment terms to regularize a standard finite difference scheme. To describe it in detail, we first fix some notations. Assume that $X$ is a $d$-dimensional hypercube, i.e. $X = \prod_{i=1}^d (a_{i},b_{i})$. We distribute $\mathcal{N}_{i}$ grid points uniformly on the $i^{\rm th}$ dimension and define
\[
h_{i} = \frac{b_{i}-a_{i}}{{\mathcal N}_{i}-1},\,\,  \mathcal{N} = \prod_{i = 1}^{d}\mathcal{N}_{i}, \mbox{ and } \Theta=\left\{\theta = (\theta_{1},\theta_{2},\cdots ,\theta_{d}){\big |}1 \le \theta_{i} \le \mathcal{N}_{i},i = 1,2,\cdots, d\right\}.
\]
Then we denote the rectangular mesh by $\mathcal{T}_{h} = \{\prod_{i=1}^d (a_{i}+(\theta_{i}-1)h_{1},a_{i}+\theta_{i}h_{1})|\theta\in\Theta\}$ and the grid points set by $X^\N = \{\pmb{x}_{\theta} = (a_{1}+(\theta_{1}-1)h_{1},a_{2}+(\theta_{2}-1)h_{2},\cdots,a_{d}+(\theta_{d}-1)h_{d})|\theta\in \Theta\}$.
the finite difference approximation of $u(\pmb{x};\pmb{\mu})$ on the grid $\mathcal{T}_{h}$ is denoted by $u^{\mathcal{N}}(\pmb{x}_{\theta};\pmb{\mu})$ represents. With appropriate rearrangement, the approximation $u^{\mathcal{N}}(\pmb{\mu})=u^{\mathcal{N}}(X^{\mathcal{N}};\pmb{\mu})$ can be regarded as a $\mathcal{N} \times 1$ vector and representing our numerical solution.

\subsubsection{Difference operators}
We first introduce several difference operators for approximating the first and second derivatives. Let $\delta_{x_{i}}^{\pm}u^{\mathcal{N}}(\pmb{x}_{\theta};\pmb{\mu})$ denote the standard forward and backward finite difference operators. That is
\[
\delta_{x_{i}}^{+}u^{\mathcal{N}}(\pmb{x}_{\theta};\pmb{\mu}):=\frac{u^{\mathcal{N}}(\pmb{x}_{\theta}+h_{i}\pmb{e}_{i};\pmb{\mu})-u^{\mathcal{N}}(\pmb{x}_{\theta};\pmb{\mu})}{h_{i}},\
\delta_{x_{i}}^{-}u^{\mathcal{N}}(\pmb{x}_{\theta};\pmb{\mu}):=\frac{u^{\mathcal{N}}(\pmb{x}_{\theta};\pmb{\mu})-u^{\mathcal{N}}(\pmb{x}_{\theta}-h_{i}\pmb{e}_{i};\pmb{\mu})}{h_{i}},
\]
where $\{\pmb{e}_{i}\}_{i=1}^{d}$ denote the canonical basis vectors for $\mathbb{R}^{d}$. The central difference operator is $\delta_{x_{i}}:=\frac{1}{2}(\delta_{x_{i}}^{+}+\delta_{x_{i}}^{-})$.
Naturally, we can define the gradient operators $\nabla_{h}^{+},\nabla_{h}^{-}$ and $\nabla_{h}$ by
\[
\nabla_{h}^{\pm}:=[\delta_{x_{1}}^{\pm},\delta_{x_{2}}^{\pm},\cdots,\delta_{x_{d}}^{\pm} ]^{T},\ \nabla_{h}:=[\delta_{x_{1}},\delta_{x_{2}},\cdots,\delta_{x_{d}} ]^{T}.
\]

Various compositions of these operators will approximate the second derivatives $\partial_{x_{i}x_{j}}^{2}$. Indeed, for $\mu,\nu\in\{+,-\}$,
$D_{h,ij}^{\mu\nu}:=\delta_{x_{j}}^{\nu}\delta_{x_{i}}^{\mu}$
are approximations of $\partial_{x_{i}x_{j}}^{2}$.
They naturally induce the approximations of the Hessian operator
\[
D_{h}^{\mu\nu}:=[D_{h,ij}^{\mu\nu}]_{i,j=1}^{d},\ \mu,\nu\in\{+,-\}.
\]
For our purpose, we will adopt the following two second order accurate approximations.
\[
\overline{D}_{h}^{2}:= \frac{D_{h}^{+-}+D_{h}^{-+}}{2},\ \widetilde{D}_{h}^{2}:= \frac{D_{h}^{++}+D_{h}^{--}}{2}.
\]

\subsubsection{The narrow-stencil finite difference scheme}
We are now ready to describe the narrow-stencil finite difference scheme for the Monge-Amp$\grave{\rm e}$re equation. We start by rewriting the Monge-Amp$\grave{\rm e}$re operator in \eqref{eq_subproblem_a} as the following form:
\[
 G(D^{2}u(\pmb{x};\pmb{\mu}),\nabla u(\pmb{x};\pmb{\mu}),u(\pmb{x};\pmb{\mu}),\pmb{x},\sigma(\pmb{\mu})) \coloneqq \sigma(\pmb{\mu})\frac{f_{X}(\pmb{x};\pmb{\mu})}{\widehat{f_{Y}}(\nabla u(\pmb{x};\pmb{\mu});\pmb{\mu})}-{\rm det}(D^{2}u(\pmb{x};\pmb{\mu}))
\]
 Then the narrow-stencil finite difference scheme seeks a grid function $u^{\mathcal{N}}(\pmb{x}_{\theta};\pmb{\mu})$ for all $\pmb{x}_{\theta}\in X^{\mathcal{N}}$ such that
$ \widehat{G}(u^{\mathcal{N}}(\pmb{x}_{\theta};\pmb{\mu}),\pmb{x}_{\theta},\sigma(\pmb{\mu})) = 0$
 where the numerical Monge-Amp$\grave{\rm e}$re operator is defined as
\begin{equation}\label{eq_numerical_operator}
\begin{aligned}
\widehat{G}(u^{\mathcal{N}}&(\pmb{x}_{\theta};\pmb{\mu}), \pmb{x}_{\theta},\sigma(\pmb{\mu}))= G(\overline{D}_{h}^{2}u^{\mathcal{N}}(\pmb{x}_{\theta};\pmb{\mu}),\nabla_{h}u^{\mathcal{N}}(\pmb{x}_{\theta};\pmb{\mu}),u^{\mathcal{N}}(\pmb{x}_{\theta};\pmb{\mu}),\pmb{x}_{\theta},\sigma(\pmb{\mu}))\\
&+2A(u^{\mathcal{N}}) : (\widetilde{D}_{h}^{2}u^{\mathcal{N}}(\pmb{x}_{\theta};\pmb{\mu})-\overline{D}_{h}^{2}u^{\mathcal{N}}(\pmb{x}_{\theta};\pmb{\mu})) -\pmb{b}(u^{\mathcal{N}})\cdot(\nabla_{h}^{+}u^{\mathcal{N}}(\pmb{x}_{\theta};\pmb{\mu})-\nabla_{h}^{-}u^{\mathcal{N}}(\pmb{x}_{\theta};\pmb{\mu})).
\end{aligned}
\end{equation}
Here $A:\mathbb{R}^{\mathcal{N}}\times X^{\mathcal{N}}\to \mathbb{R}^{d\times d}$ is a matrix-valued function and $\pmb{b}:\mathbb{R}^{\mathcal{N}}\times X^{\mathcal{N}}\to \mathbb{R}^{d}$ vector-valued. In this article, we simply choose $A(u^{\mathcal{N}}(\pmb{x}_{\theta};\pmb{\mu}),\pmb{x}_{\theta}) = \alpha I$ and $\pmb{b}(u^{\mathcal{N}}(\pmb{x}_{\theta};\pmb{\mu}),\pmb{x}_{\theta}) = \beta \pmb{e}$ for $\alpha \ge0$ and $\beta \ge0$, where $I$ denotes the $d \times d$ identity matrix and $\pmb{e}$ the $d$-dimensional column vector with all elements equal to $1$.

Next we introduce the discretization of the boundary condition and uniqueness condition \eqref{eq_subproblem_b}  and \eqref{eq_subproblem_c}. Since (\ref{eq_numerical_operator}) is an approximation of the fourth order PDE, we have to introduce an additional boundary condition to guarantee that the discrete problem is well-posed. As \cite{feng2009vanishing} has done, here we introduce one discrete additional boundary condition
\begin{equation}\label{BC_Neumann_additional_discretization}
\nabla_{h}(\Delta_{h} u^{\mathcal{N}}(\pmb{x}_{\theta};\pmb{\mu}))\cdot \pmb{n}(\pmb{x}_\theta) = 0,\   \pmb{x}_{\theta}\in X^{\mathcal{N}}\cap\partial X,
\end{equation}
where discrete operator $\Delta_{h}$ is defined by the nine points finite difference scheme. For the Neumann boundary condition, we simply discretize it by the central difference
\begin{equation}\label{BC_Neumann_discretization}
\nabla_{h}u^{\mathcal{N}}(\pmb{x}_{\theta};\pmb{\mu})\cdot\pmb{n}(\pmb{x}_\theta) = \phi(\pmb{x}_{\theta}),\   \pmb{x}_{\theta}\in X^{\mathcal{N}}\cap\partial X.
\end{equation}
Finally, the uniqueness condition is simply approximated by the mean value of each component of vector $u^{\mathcal{N}}(\pmb{\mu})$, which is defined as
\begin{equation}\label{discrete_uniquenenss_condition}
\frac{\sum_{\pmb{x}_{\theta}\in X^{\mathcal{N}}} u^{\mathcal{N}}(\pmb{x}_{\theta};\pmb{\mu})}{\mathcal{N}}=0.
\end{equation}

The above procedure presents the discretization of the subproblem (\ref{eq_subproblem_a}-\ref{eq_subproblem_c}), by which we can obtain a system of equations. The equations can be solved by the Newton's method efficiently, with details on initialization provided in Section \ref{Section_The algorithm for transport problem}.

\begin{wrapfigure}[15]{R}{0.45\textwidth}
\centering
\includegraphics[scale=0.39]{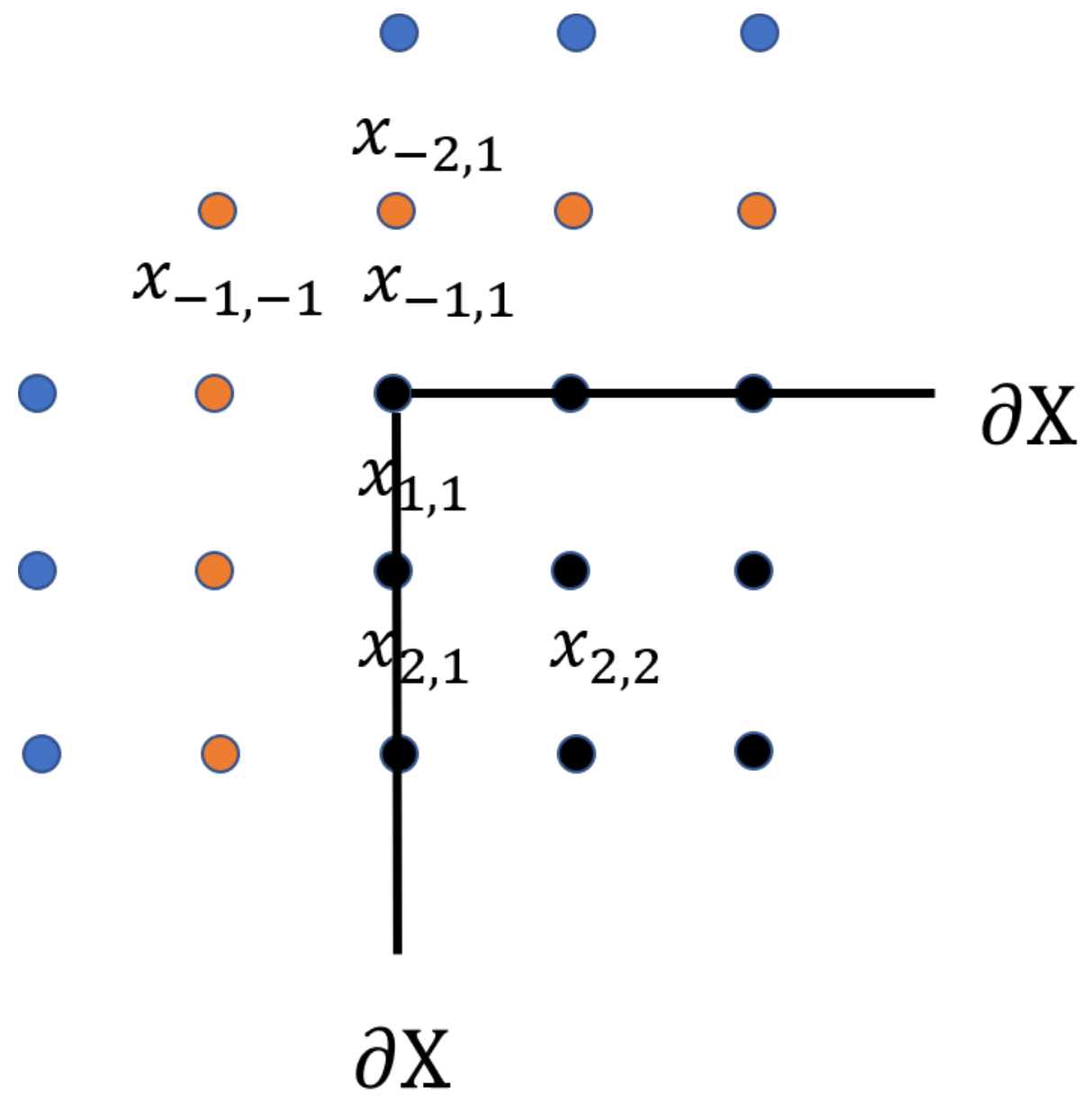}
\caption{Schematic of the ghost points.}
\label{ghost_points}
\end{wrapfigure}
We note here that handling the boundary conditions (\ref{BC_Neumann_additional_discretization}, \ref{BC_Neumann_discretization}) will require the introduction of  two layers of ghost points as depicted in Figure \ref{ghost_points}. The values at the first layer of ghost points near the boundary are determined by the Neumann boundary condition. For instance, the normal derivative in the direction $\pmb{n}_{x_{2}}$ at point $\pmb{x}_{1,1}$ can be discretized as
\[
u_{\pmb{n}_{x_{2}}}(\pmb{x}_{1,1};\pmb{\mu}) \approx \frac{1}{2h_{2}}(u^{\mathcal{N}}(\pmb{x}_{-1,1};\pmb{\mu})-u^{\mathcal{N}}(\pmb{x}_{2,1};\pmb{\mu})).
\]
In the general case, we define the normal derivative in the diagonal direction at the four corners by the sum of the normal derivative in two orthogonal outward directions and still apply the the central difference scheme to discretize the derivative.
\[
u_{\pmb{n}_{\rm diag}}(\pmb{x}_{1,1};\pmb{\mu}) \approx \frac{1}{2\sqrt{h_{1}^{2}+h_{2}^{2}}}( u^{\mathcal{N}}(\pmb{x}_{-1,-1};\pmb{\mu})-u^{\mathcal{N}}(\pmb{x}_{2,2};\pmb{\mu})) =\frac{1}{\sqrt{2}}(u_{\pmb{n}_{x_{1}}}(\pmb{x}_{1,1};\pmb{\mu})+u_{\pmb{n}_{x_{2}}}(\pmb{x}_{1,1};\pmb{\mu})).
\]
The values of the most outer ghost points are determined by the additional boundary condition \eqref{BC_Neumann_additional_discretization} since for example
\[
0 = \frac{\partial\Delta_{h}u^{\mathcal{N}}(\pmb{x}_{1,1};\pmb{\mu})}{\partial\pmb{n}_{x_{2}}}
= \frac{(\frac{u^{\mathcal{N}}(\pmb{x}_{-2,1};\pmb{\mu})-2u^{\mathcal{N}}(\pmb{x}_{-1,1};\pmb{\mu})+u^{\mathcal{N}}(\pmb{x}_{1,1};\pmb{\mu})}{h_{2}^{2}}+\delta_{x_{1}}^{2}u^{\mathcal{N}}(\pmb{x}_{-1,1};\pmb{\mu}))-(\delta_{x_{1}}^{2}+\delta_{x_{2}}^{2})u^{\mathcal{N}}(\pmb{x}_{2,1};\pmb{\mu})}{2h_{2}}
\]
means that
\[
u^{\mathcal{N}}(\pmb{x}_{-2,1};\pmb{\mu}) = -u^{\mathcal{N}}(\pmb{x}_{1,1};\pmb{\mu})+2u^{\mathcal{N}}(\pmb{x}_{-1,1};\pmb{\mu})+h_{2}^{2}(-\delta_{x_{1}}^{2}u^{\mathcal{N}}(\pmb{x}_{-1,1};\pmb{\mu})+(\delta_{x_{1}}^{2}+\delta_{x_{2}}^{2})u^{\mathcal{N}}(\pmb{x}_{2,1};\pmb{\mu})).
\]

\begin{remark}
This scheme only entails a 14-point stencil in two dimension. By contrast, the monotone FDM of \cite{froese2011convergent} needs more points because of the discrete-direction error. For this reason, it is called the narrow-stencil scheme. We next examine the last two terms in the scheme \eqref{eq_numerical_operator} to reveal that the narrow-stencil FDM introduces the stabilization terms, i.e., numerical moment and numerical viscosity \cite{feng2009mixed,feng2011vanishing,1984Two}. Indeed, a direct calculation shows that
\[
(\delta_{x_{i}}^{+}-\delta_{x_{i}}^{-})u^{\mathcal{N}}(\pmb{x}_{\theta};\pmb{\mu}) = h_{i}\frac{u^{\mathcal{N}}(\pmb{x}_{\theta}-h_{i}\pmb{e}_{i};\pmb{\mu})-2u^{\mathcal{N}}(\pmb{x}_{\theta};\pmb{\mu})+u^{\mathcal{N}}(\pmb{x}_{\theta}+h_{i}\pmb{e}_{i};\pmb{\mu})}{h_{i}^{2}}=h_{i}\delta_{x_{i}}^{2}u^{\mathcal{N}}(\pmb{x}_{\theta};\pmb{\mu}).
\]
Therefore
$\beta\pmb{e}\cdot(\nabla_{h}^{+}-\nabla_{h}^{-})u^{\mathcal{N}}(\pmb{x}_{\theta};\pmb{\mu})\approx \beta h\Delta u(\pmb{x}_{\theta}\pmb{\mu})$,
which amounts to addition of numerical viscosity, a known technique for constructing convergent difference scheme, see e.g. \cite{1984Two}. Further, we can show that, for all $i,j\in \{1,2,\cdots, d\}$
\[
(\widetilde{D}_{h,ij}^{2}-\overline{D}_{h,ij}^{2})u^{\mathcal{N}}(\pmb{x}_{\theta};\pmb{\mu}) = \frac{1}{2}(\delta_{x_{i}}^{+}\delta_{x_{j}}^{+}+\delta_{x_{i}}^{-}\delta_{x_{j}}^{-}-\delta_{x_{i}}^{+}\delta_{x_{j}}^{-}-\delta_{x_{i}}^{-}\delta_{x_{j}}^{+})u^{\mathcal{N}}(\pmb{x}_{\theta};\pmb{\mu}) =\frac{h_{i}h_{j}}{2}\delta_{x_{i}}^{2}\delta_{x_{j}}^{2}u^{\mathcal{N}}(\pmb{x}_{\theta};\pmb{\mu}).
\]
This means that
\[
\alpha I:(\widetilde{D}_{h}^{2}-\overline{D}_{h}^{2})u^{\mathcal{N}}(\pmb{x}_{\theta};\pmb{\mu})\approx \alpha h^{2}\Delta^{2} u(\pmb{x}_{\theta};\pmb{\mu}),
\]
which introduces the numerical moment in the vanishing moment method of \cite{feng2009mixed}. Synthesizing the above observations, it is clear that the proposed scheme is an approximation of the following fourth order quasilinear PDE
\[
\alpha h^{2} \Delta^{2}u(\pmb{x};\pmb{\mu}) -\beta h \Delta u(\pmb{x};\pmb{\mu})+G(D^{2}u(\pmb{x};\pmb{\mu}),\nabla u(\pmb{x};\pmb{\mu}),u(\pmb{x};\pmb{\mu}),\pmb{x},\sigma(\pmb{\mu})) = 0,
\]
which is a regularization of the original nonlinear PDE.
\end{remark}

\subsection{The algorithm for the transport boundary problem }\label{Section_The algorithm for transport problem}
Integrating the iterative approach for the transport boundary condition with the narrow-stencil finite difference scheme for the Neumann boundary subproblem, we are ready to present our full order model for the transport boundary problem of the Monge-Amp$\grave{\rm e}$re equation. For notational simplicity, we define the nonlinear approximate equations (\ref{eq_numerical_operator}) with the boundary condition (\ref{BC_Neumann_discretization}), (\ref{BC_Neumann_additional_discretization})
and the uniqueness condition (\ref{discrete_uniquenenss_condition}) as the following system of equations
\begin{equation}\label{eq_subproblem_discrete}
F(u^{\mathcal{N}}(\pmb{\mu}),\sigma(\pmb{\mu});\phi(X^{\mathcal{N}}\cap \partial X;\pmb{\mu})) = 0,
\end{equation}
where $F(\cdot,\cdot;\cdot))$ is a nonlinear system with $(\mathcal{N}+1)$ equations on $(u^{\mathcal{N}}(X^{\mathcal{N}};\pmb{\mu}),\sigma(\pmb{\mu}))^T$. At each boundary iteration in Section \ref{Section:Transport boundary condition}, the Newton's method is used to solve the above system of equations which is recalled as
\[
F(u^{\mathcal{N},k}(\pmb{\mu}),\sigma(\pmb{\mu});\phi^k(X^{\mathcal{N}}\cap \partial X;\pmb{\mu})) = 0.
\]
To assist with the convergence of the Newton's method, an initialization sufficiently close to the exact solution is necessary.
We adopt the approach of \cite{froese2011convergent} and take the initial value of the $k$-th iteration,  $(u_0^{\mathcal{N},k}(X^{\mathcal{N}};\pmb{\mu}),\sigma_0^k(\pmb{\mu}))^T$, as the solution of
\[
\begin{aligned}
&\Delta_h u^\N(\pmb{x}_\theta;\pmb{\mu}) = \sigma(\pmb{\mu})\left(\frac{2f_{X}(\pmb{x}_\theta;\pmb{\mu})}{f_{Y}(\pmb{x}_\theta-\pmb{y}_{0};\pmb{\mu})}\right)^{\frac{1}{2}},\ \pmb{x}_\theta\in X^\N,\\
&\frac{\sum_{\pmb{x}_{\theta}\in X^{\mathcal{N}}} u^{\mathcal{N}}(\pmb{x}_{\theta};\pmb{\mu})}{\mathcal{N}}=0,\\
&\nabla_{h}u^{\mathcal{N}}(\pmb{x}_{\theta};\pmb{\mu})\cdot\pmb{n}(\pmb{x}_\theta) = \phi^k(\pmb{x}_{\theta}),\   \pmb{x}_{\theta}\in X^{\mathcal{N}}\cap\partial X,
\end{aligned}
\]
We are now ready to present the algorithm for solving the problem \eqref{eq_MA} for a given value of $\pmb{\mu}$, in Algorithm \ref{alg:1}.
\begin{algorithm}[H]
\caption{The algorithm for the transport boundary problem of Monge-Amp$\grave{\rm e}$re equation}
\label{alg:1}
\begin{algorithmic}[1]
 \STATE Set $\phi^{-1}(\pmb{x};\pmb{\mu}) = B\pmb{x}\cdot\pmb{n}$ for sufficient large $B$, the error tolerance $\epsilon$ and the maximum number of iterations $K$. Compute $(u^{\mathcal{N},0}(\pmb{\mu}), \sigma^{0}(\pmb{\mu}))^T$ by solving $F(u^{\mathcal{N}}(\pmb{\mu}),\sigma(\pmb{\mu});\phi^{-1}(X^{\mathcal{N}}\cap \partial X;\pmb{\mu})) = 0$.
 \STATE Initialize $k = 0$ and $r = 1$.
 \WHILE{$r\ge \epsilon$ {and} $k<K$}
 \STATE Compute $\phi^{k}(\pmb{x}_{\theta};\pmb{\mu})= \mathbb{P}_{\partial Y}(\nabla_{h} u^{\mathcal{N},k}(\pmb{x}_{\theta};\pmb{\mu}))\cdot \pmb{n}$, for $\pmb{x}_{\theta}\in X^{\mathcal{N}}\cap\partial X $.
 \STATE Compute $(u^{\mathcal{N},k+1}(\pmb{\mu}), \sigma^{k+1}(\pmb{\mu}))^T$ by solving $F(u^{\mathcal{N}}(\pmb{\mu}),\sigma(\pmb{\mu});\phi^{k}(X^{\mathcal{N}}\cap \partial X;\pmb{\mu})) = 0$.
 \STATE Compute the relative error $r = ||u^{\mathcal{N},k+1}(\pmb{\mu})-u^{\mathcal{N},k}(\pmb{\mu})||_{\ell^{\infty}(\mathbb{R}^{\mathcal{N}})}$.
 \STATE Let $k = k+1$.
 \ENDWHILE
\end{algorithmic}
\end{algorithm}

\section{The Reduced Residual Reduced Over-Collocation (R2-ROC) Method}
\label{sec:R2-ROC-Alg}

Following the full order model presented in the last section, we introduce our proposed Reduced Order Model (ROM) for the transport boundary problem of the Monge-Amp$\grave{\rm e}$re equation. Specifically, we adopt the Reduced Over Collocation (ROC) approach \cite{chen2021eim} developed for  parametrized nonlinear partial differential equations. The unique feature is the immunity of the degradation in online efficiency suffered by classical RBM as a result of the EIM-like expansion of the nonlinear and nonaffine terms. {To illustrate the algorithm, we first recall that we denote by $u(\pmb{x}; \bmu)$ the exact solution of the Monge-Amp$\grave{\rm e}$re equation (\ref{eq_MA}) which is nonlinear and parameterized in a nonaffine fashion by $\pmb{\mu}$.
Moreover, the resulting FOM solution corresponding to parameter $\bmu$ is denoted by $u^\N(X^{\N}; \bmu)$ which we assume is close enough to the exact solution $u(\pmb{x}; \bmu)$ for us to adopt as a reference for the ROM.  }
Now we are ready to briefly review the R2-ROC algorithm. It has two components: an online (reduced) solver of size $n$ that is between $1$ and $N$ with $N$ usually much smaller than $\N$, and an offline training component which repeatedly calls the online solver of increasing size $n$ to build up a surrogate solution space from scratch dimension-by-dimension.

\subsubsection*{Online solver}

Given the reduced space $W_n$ and a collocation set $X^m$, a subset of the full grid $X^\N$, of cardinality $m$ that is comparable to $n$, R2-ROC identifies a surrogate solution for any specific parameter $\bmu$ in the following form
\[
\widehat{u}_n(\bmu) = W_{n} \bc_n (\bmu).
\]
Here, for simplicity of notation, we also adopt $W_{n}$ for the snapshot matrix whose column space forms the reduced space $W_{n}$.
We then subject this surrogate solution to the FOM equation \eqref{eq_subproblem_discrete} which encompasses  equations \eqref{eq_numerical_operator} - \eqref{discrete_uniquenenss_condition}. Note however that $\widehat{u}_n(\bmu)$ automatically satisfies the uniqueness condition \eqref{discrete_uniquenenss_condition} due to that the constraint
\[
\frac{\sum_{\pmb{x}_{\theta}\in X^{\mathcal{N}}} \widehat{u}_{n}(\pmb{x}_{\theta};\pmb{\mu})}{\mathcal{N}}\equiv0.
\]
is linear and that all RB snapshots $\{u_{i}\}_{i=1}^{n}$ satisfy it by definition. Therefore, we just need to subject $\widehat{u}_{n}(\pmb{\mu})$ to \eqref{eq_numerical_operator} - \eqref{BC_Neumann_discretization} a nonlinear system
\[
F_{r}(W_n \pmb{c}_{n}(\pmb{\mu}),\sigma(\pmb{\mu});\phi^k(X^{\mathcal{N}}\cap \partial X;\pmb{\mu})) = 0,
\]
with $\mathcal{N}$ equations for unknown $(\pmb{c}_{n}(\pmb{\mu}), \sigma(\pmb{\mu}))^T$. The $k^{\rm th}$ iterate of the solution, $\bc_n^k(\bmu) \in {\mathbb R}^{n \times 1}$ and $\sigma^k(\pmb{\mu})$, is obtained by minimizing a subsampled residual
\[
(\bc_n^k(\bmu),\sigma^k(\pmb{\mu}))^T = \argmin \left\lVert P_* \left(F_{r}(W_n \pmb{c}_{n}(\pmb{\mu}),\sigma(\pmb{\mu});\phi^k(X^{\mathcal{N}}\cap \partial X;\pmb{\mu})) \right)\right\rVert_{\ell^2(\mathbb{R}^m)}.
\]
The RB space $W_{n}$,  the reduced collocation set $X^m$, and subsampling matrix $P_* \in \mathbb{R}^{m \times \N}$ that is constructed according to $X^m$, will be generated in the offline process that is described next.
The online algorithm is presented in Algorithm \ref{alg:2}.

\begin{algorithm}[htb]
\caption{Online algorithm: The reduced algorithm for the transport boundary value problem of Monge-Amp$\grave{\rm e}$re equation}
\label{alg:2}
\begin{algorithmic}[1]
 \STATE Set $\phi^{-1}(\pmb{x};\pmb{\mu}) = B\pmb{x}\cdot\pmb{n}$ for sufficient large $B$, the error tolerance $\epsilon$ and the maximum number of iterations $K$. Compute $(\pmb{c}_{n}^{0}(\pmb{\mu}), \sigma^{0}(\pmb{\mu}))^T$ by solving $P_{\star}F_1(W_n \pmb{c}_{n}(\pmb{\mu}),\sigma(\pmb{\mu});\phi^{-1}(X^{\mathcal{N}}\cap \partial X;\pmb{\mu})) = 0$.
 \STATE Initialize $k = 0$ and $r = 1$.
 \WHILE{$r\ge\epsilon\, {\mbox{ and }} \,k< K$}
 \STATE Compute $\phi^{k}(\pmb{x}_{\theta};\pmb{\mu})= \mathbb{P}_{\partial Y}(\nabla_{h} \widehat{u}_{n}^{k}(\pmb{\mu})(\pmb{x}_{\theta};\pmb{\mu}))\cdot \pmb{n}$, for some $\pmb{x}_{\theta}\in S^k(\bmu)$, where $\widehat{u}_{n}^{k}(\pmb{\mu}) = W_n \pmb{c}_{n}^{k}(\pmb{\mu})$.\label{iteration_BC}
 \STATE Compute $(\pmb{c}_{n}^{k+1}(\pmb{\mu}), \sigma^{k+1}(\pmb{\mu}))^T$ by solving $\left\lVert P_{\star}F_1(W_n \pmb{c}_{n}(\pmb{\mu}),\sigma(\pmb{\mu});\phi^{k}(X^{\mathcal{N}}\cap \partial X;\pmb{\mu}))\right\rVert_{\ell^2(\mathbb{R}^m)} = 0$.\label{c}
 \STATE Compute the relative error $r = \left\lVert\pmb{c}_n^{k+1}(\pmb{\mu})-\pmb{c}_n^{k}(\pmb{\mu})\right\rVert_{\ell^{\infty}(\mathbb{R}^{n})}$.
 \STATE Let $k = k+1$.
 \ENDWHILE
\end{algorithmic}
\end{algorithm}

\subsubsection*{Online efficiency and robustness with respect to the shortest-distance projection $\mathbb{P}_{\partial Y}$}
The RB method is said to be online-efficient if the RB solver can be assembled and RB approximation solved in complexity independent of $\mathcal{N}$ in the online stage and the error estimator can be computed, via an offline-online decomposition if necessary, in complexity independent of $\mathcal{N}$ online \cite{casenave2014accurate}.
R2-ROC method is online-efficient as established by \cite{chen2021eim}. Our version of R2-ROC for solving the parameterized Monge-Amp$\grave{\rm e}$re equation with a transport boundary features the added layer of iteration and
 the shortest-distance projection ${\mathbb P}_{\partial Y}$ in \eqref{eq_Neumann_BC_of_subproblem}. We note that the iteration is up to a fixed number $K$ and that ${\mathbb P}_{\partial Y}$ is only carried out for part of the boundary points, $S^k(\bmu)\, (\subset X^{\mathcal{N}}\cap\partial X)$, whose cardinality only depends on $m$. Therefore, we conclude that our R2-ROC remains online-efficient.

Moreover, we emphasize that the shortest-distance projection $\mathbb{P}_{\partial Y}(\nabla u^{k}(\pmb{x};\pmb{\mu}))$ at each iteration must be calculated whenever possible exactly, that is, without discretizing $\partial Y$. This is straightforward when, for example, the target boundary $\partial Y$ is a polygon or a circle. In fact if $\partial Y$ is discretized as $\partial Y_{\N}$, that is we calculate $\mathbb{P}_{\partial Y}(\pmb{v})$ as $\argmin_{\pmb{y} \in \partial Y_{\N}} \lVert \pmb{y} - \pmb{v}\rVert$, our reduced solver is much less robust.

\subsubsection*{Offline training}

The offline component utilizes the classical parameter-greedy framework with an error indicator based on the hyper-reduced residual \cite{chen2021eim} to iteratively construct  the reduced basis space $W_n$ and subsequently enrich the collocation set $X^m$ which determines the subsampling matrix $P_*$. The algorithm judiciously identifies parameter values
$\left\{\bmu^1, \dots, \bmu^N \right\}$
one-by-one and construct the reduced basis space via the corresponding snapshots.
With these notations set, we start the greedy procedure with a randomly chosen $\bmu^1$ and obtain the
snapshots $u^\N(X^\N; \bmu^1)$ by the high fidelity algorithm.
The RB space $W_1$ is then set as $W_1 =\{u_1\}= \{u^\N(X^\N; \bmu^1)\}$, and the first collocation point chosen as the EIM point of the first basis $ \bx_\ast^1 = \argmax_{x \in X^\N} |u_1|$. We then use the online procedure described above to obtain an RB approximation $\widehat{u}_{n}(\boldsymbol{\mu})$ for each parameter $\bmu$ in $\Xi_{\rm train}$ (a discretization of the parameter domain $\calD$)
and compute its error estimator $\Delta_n^{RR}(\bmu)$.
\[
\Delta_n^{RR}(\bmu) \coloneqq  \lVert P_\ast r_n(\bmu)\rVert_{\ell^\infty}.
\]
Here, $r_n(\bmu) = F_1(W_n \pmb{c}_{n}^{K_n(\bmu)}(\pmb{\mu}),\sigma^{K_n(\bmu)}(\pmb{\mu});\phi^{K_n(\bmu)-1}(X^{\mathcal{N}}\cap \partial X;\pmb{\mu}))$ is the {\em full} residual for the current RB approximation $\widehat{u}_n(\bmu)$ of parameter $\bmu$, and $K_n(\bmu)$ is the corresponding number of iterations.
$P_\ast r_n(\bmu)\in \mathbb{R}^{m \times 1}$ then represents its {\em reduced} (subsampled) version \footnote{The conventional error estimate calculates the negative-order norm of the residual and scales it by the (parametric) stability factor. It is challenging to compute for the nonlinear and nonaffine case with EIM expansion due to the involvement of the successive constraint method \cite{huynh2007successive,huynh2010natural} used to efficiently estimate the parametric stability factor, and the delicacy of evaluating the residual norm even for the linear problem \cite{casenave2014accurate,chen2019robust}. This simple error estimator based on the reduced residual \cite{chen2021eim} has shown to be promising for nonlinear and nonaffine problems without the need of EIM expansion.} whose evaluation is independent of $\calN$.
After these error indicators are evaluated, we proceed as follows to enrich the RB space and expand the collocation sets.
\begin{enumerate}
\item {\bf Greedy in $\bmu$:}
The greedy choice  is through maximizing $\Delta_n^{RR}(\bmu)$ over the training set $\Xi _{\rm train}$:
\[
  \bmu^{n+1} = \argmax_{\bmu \in \Xi _{\rm train}} \Delta_n^{RR}(\bmu). 
\]

\item {\bf $X^m$ expansion:} With the newly selected $\bmu^{n+1}$, we solve for the truth approximations $u_{n+1}$. We then obtain the first additional collocation point from the EIM process of $u_{n+1}$, and the second additional point by the EIM process of the full residual $r_n(\bmu^{n+1})$.

\end{enumerate}
The offline algorithm of the ROC method is shown in Algorithm \ref{alg:3}. We refer the readers to \cite{chen2021eim} for more details including the analysis which demonstrates the importance of retaining these two sets of points for producing accurate approximations for both the reduced solution and the residual corresponding to each parameter.

\begin{algorithm}[htb]
\caption{Offline algorithm: the reduced over-collocation methods the transport boundary value problem of Monge-Amp$\grave{\rm e}$re equation}
\label{alg:3}
\begin{algorithmic}[1]
\STATE Choose $\pmb{\mu}^{1}$ randomly in $\Xi_{\rm train}$ and obtain $u^{\mathcal{N}}(\pmb{\mu}^{1})$ by Algorithm \ref{alg:1}. Find $\pmb{x}_{\star}^{2} = \mathop{\arg\max}_{\pmb{x}_{\theta}\in X^{\mathcal{N}}}|u^{\mathcal{N}}(\pmb{x}_{\theta};\pmb{\mu}^{1})|$ and
$\pmb{x}_{\star}^{1} = \mathop{\arg\max}_{\pmb{x}_{\theta}\in X^{\mathcal{N}}}|u^{\mathcal{N}}(\pmb{x}_{\theta};\pmb{\mu}^{1})-u^{\mathcal{N}}(\pmb{x}_{\star}^{2};\pmb{\mu}^{1})|$. Then let $n = 1$, $X^{m} = X_{s}^{n} = \{\pmb{x}_{\star}^{1},\pmb{x}_{\star}^{2}\}$, and $u_{1} = u^{\mathcal{N}}(\pmb{\mu}^{1})/u^{\mathcal{N}}(\pmb{x}_{\star}^{2};\pmb{\mu}^{1})$. \STATE Initialize $W_{1} = \{u_{1}\}$ and $X_{r}^{0} = \emptyset $.
\FOR{$n = 2,\cdots, N$}
\STATE Solve $\pmb{c}_{n-1}(\pmb{\mu})$ by Algorithm \ref{alg:2} with $W_{n-1}, X^{m}$ and obtain the corresponding number of iterations $K_n(\bmu)$ and calculate $\Delta_{n-1}^{RR}$ for every $\pmb{\mu}\in \Xi_{\rm train}$.\label{alg1_RB_approximation}
\STATE Find $\pmb{\mu}^{n} = \mathop{\arg\max}_{\pmb{\mu}\in \Xi_{\rm train}}\Delta_{n-1}(\pmb{\mu})$.
\STATE Solve $u_{n}:=u^{\mathcal{N}}(\pmb{\mu}^{n})$ by Algorithm \ref{alg:1}. Orthogonalize $u_{n}$: find $\{\alpha_{j}\}$ and let $u_{n} = u_{n}-\sum_{j=1}^{n-1}\alpha_{j}u_{j}$ such that $u_{n}(X_{s}^{n-1}\backslash\{\pmb{x}_\star^1\})=0$.
\STATE Find $\pmb{x}_{\star}^{n+1} = \mathop{\arg\max}_{\pmb{x}_{\theta}\in X^{\mathcal{N}}}|u_{n}(\pmb{x}_{\theta})|$, $u_{n} = u_{n}/u_{n}(\pmb{x}_{\star}^{n+1})$, and let $X_{s}^{n} = X_{s}^{n-1}\cup\{\pmb{x}_{\star}^{n+1}\}$.
\STATE Assume that $\sigma(\pmb{\mu}^{n}), \phi(X^{\mathcal{N}}\cap \partial X;\pmb{\mu}^{n})$ are obtained when solving $\pmb{c}_{n-1}(\pmb{\mu}^{n})$ in Step \ref{c} of Algorithm \ref{alg:2}. Then compute the full residual vector $r_{n-1} = F_1(W_{n-1} \pmb{c}_{n-1}(\pmb{\mu}^{n}),\sigma(\pmb{\mu}^{n});\phi(X^{\mathcal{N}}\cap \partial X;\pmb{\mu}^{n}))$ and orthogonalize $r_{n-1}$: find $\{\alpha_{j}\}$ and let $r_{n-1} = r_{n-1}-\sum_{j=1}^{n-2}\alpha_{j}r_j$ so that $r_{n-1}(X_{r}^{n-2})=0$.
\STATE Find $\pmb{x}_{\star\star}^{n-1} = \mathop{\arg\max}_{\pmb{x}_{\theta}\in \pmb{x}^{\mathcal{N}}}|r_{n-1}(\pmb{x}_{\theta})|$. Let $r_{n-1} = r_{n-1}/r_{n-1}(\pmb{x}_{\star\star}^{n-1})$ and $X_{r}^{n-1} = X_{r}^{n-2}\cup\{\pmb{x}_{\star\star}^{n-1}\}$.
\STATE Update $W_{n} = [W_{n-1},u_{n}]$, $X^{m} = X_{s}^{n}\cup X_{r}^{n-1}$.
\ENDFOR
\end{algorithmic}
\end{algorithm}

\section{Numerical results}
\label{Section:numerical}
In this section, we present the computational results to verify two main works of this paper: the proposed narrow-stencil FDM is effective in solving the transport boundary problem of the Monge-Amp$\grave{\rm e}$re equation and the adapted R2-ROC method can efficiently solve the parameterized transport boundary problem.

\subsection{FDM results}
In this section, we focus on our Monge-Amp$\grave{\rm e}$re equation solver. The results of four tests are presented to gauge the performance of the proposed approach for approximating the viscosity solutions. The problems are described in Table \ref{tab:fdm_test_t} which lists the original density $f_{X}(\pmb{x})$, the target density $f_{Y}(\pmb{y})$, and the exact mapping $\nabla_{\pmb{x}} u$. The first one maps the square $(-0.5,0.5)\times (-0.5,0.5)$ onto the rectangle$(0.5,1.5)\times (-0.25,0.25)$. The second one maps the square $(-0.5,0.5)\times(-0.5,0.5)$ onto the same square where we define the following auxiliary function
\[
q(z) = (-\frac{1}{8\pi}z^2+\frac{1}{256\pi ^{3}}+\frac{1}{32\pi}){\cos}(8\pi z)+\frac{1}{32\pi^{2}}z{\sin}(8\pi z).
\]
The third one maps a uniform density on the unit square $(0,1)\times(0,1)$ onto a density that blows up at a point on the same square. The last one maps a uniform density on the square $(-0.5,0.5)\times (-0.5,0.5)$ onto a Gaussian density on the disk $y_1^2+y_2^2\le 0.5^2$. This last test is meant to verify the effectiveness of our approach for transporting a rectangular  boundary to a circular one, a nontrivial task.
\begin{table}[htb]
\resizebox{\textwidth}{!}{
    \centering
    \renewcommand{\arraystretch}{1.5}
    \begin{tabular}{|c|c|c|c|}
    \hline
    Test& Original density & Target density & Exact mapping\\
    \hline
    1&$f_{X}(\pmb{x}) =\frac{1}{0.16}\exp\left(-\frac{1}{2}\frac{x_{1}^{2}}{0.4^{2}}-\frac{1}{2}\frac{x_{2}^{2}}{0.4^{2}}\right)$     & $f_{Y}(\pmb{y}) =\frac{1}{0.08}\exp\left(-\frac{1}{2}\frac{(y_{1}-1)^{2}}{0.4^{2}}-\frac{1}{2}\frac{y_{2}^{2}}{0.2^{2}}\right)$ & $\nabla_{\pmb{x}} u=\left(x_{1}+1,\frac{x_{2}}{2}\right)$\\
    \hline
\multirow{2}{*}{2}&$f_{X}(\pmb{x})=1+4(q''(x_{1})q(x_{2})+q(x_{1})q''(x_{2}))+$         & \multirow{2}{*}{$1$} & $u_{x_{1}}=x_{1}+4q'(x_{1})q(x_{2})$\\
 &$16(q(x_{1})q(x_{2})q''(x_{1})q''(x_{2})-q'(x_{1})^{2}q'(x_{2})^{2})$         && $u_{x_{2}}=x_{2}+4q'(x_{2})q(x_{1})$\\
    \hline
    3&$1$ & $f_{Y}(\pmb{y}) = \frac{\exp(-2\sqrt{(y_{1}-0.5)^{2}+(y_{2}-0.5)^{2}})}{\sqrt{(y_{1}-0.7)^{2}+(y_{2}-0.7)^{2}}}$ & -\\
    \hline
    4&1     & $f_{Y}(\pmb{y}) =1+\frac{1}{0.02\pi}\exp\left(-\frac{y_{1}^{2}+y_{2}^{2}}{0.02}\right)$ & -\\
    \hline
    \end{tabular}
}
    \caption{Setup of the test problems for the transport boundary case.}
    \label{tab:fdm_test_t}
\end{table}

\begin{figure}[htb]
\centering
\includegraphics[scale=0.5]{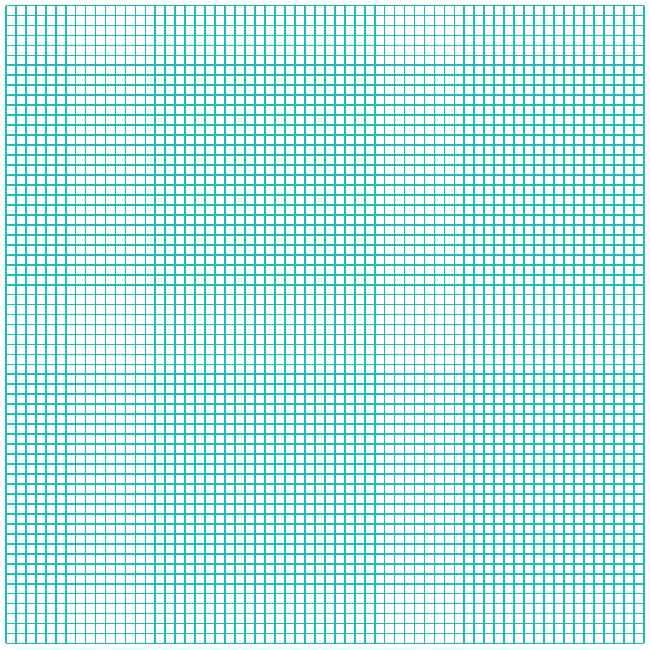}
\includegraphics[scale=0.5]{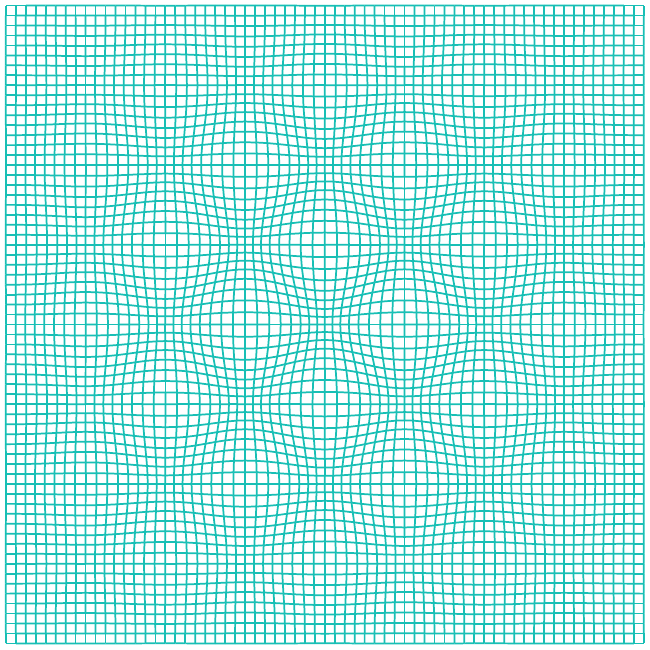}
\includegraphics[scale=0.5]{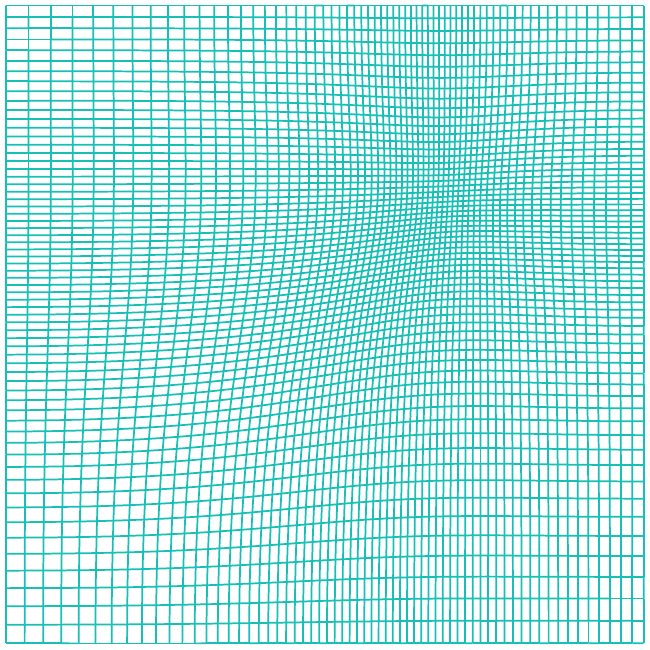}
\includegraphics[scale=0.5]{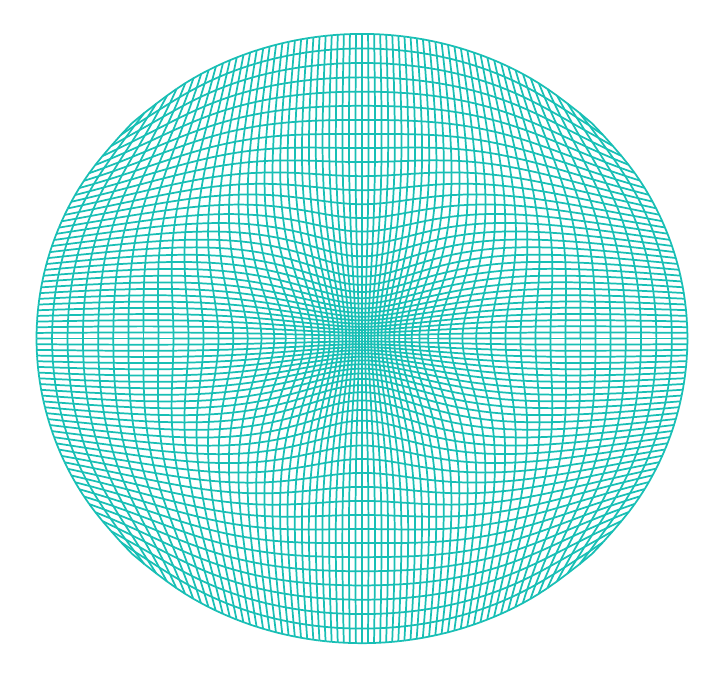}
\caption{A uniform Cartesian mesh and its images under the map $\nabla u$ for the second and third test.}
\label{fig:fdm_t_images}
\end{figure}

We present the result in Table \ref{tab:fdm_T}. {Here, the maximum errors are computed based on the exact solutions or the solutions on the finest grids.}
From the table we can see the solution indeed achieve machine accuracy for the first test, order $2$ accuracy for the second even without $\alpha$ and $\beta$. For the third, many currently available methods become slow or unstable when the ratio $R = \min\{\frac{f_{X}(\pmb{x})}{f_{Y}(\nabla u(\pmb{x}))}\}$ is small. We see that our approach works well, with small $\alpha$, even when $R$ is very small. We also provide, in Figure \ref{fig:fdm_t_images}, a uniform Cartesian mesh and its images under the second, third and fourth map.
\begin{table}[htb]
\resizebox{\textwidth}{!}{
 \centering
\bigskip
\centering
\begin{small}
\begin{tabular}{ | c |    c | c|c|}
\Xhline{1pt}
\multicolumn{1}{ | c |}{\multirow{2}{*}{$\mathcal{N}$}}
&\multicolumn{2}{ c |}{Test 1 ($\alpha = \beta = 0$)}\\ \cline{2-3}
\multicolumn{1}{ | c |}{}
&time(s) & max error	\\\hline
$15^{2}$&    0.48&   2.44E-15        \\
$31^2$&    0.11&   5.88E-15       \\
$65^2$&   0.75&   7.11E-15       \\
$127^2$&   7.71&   2.89E-14         \\
$255^2$&   90.70&   5.68E-14     \\
$511^2$&   2789.93&   7.17E-14     \\\Xhline{1pt}
\end{tabular}
\begin{tabular}{  | c|c |c|}
\Xhline{1pt}
\multicolumn{3}{| c |}{Test 2 ($\alpha = \beta = 0$)}	\\\cline{1-3}
time(s) & max error& order	\\\hline
 0.05&   2.72E-03   & --     \\
0.06&   7.47E-04   &1.78     \\
 0.51&   2.24E-04   &1.63     \\
 3.77&   6.18E-05   &1.92      \\
 51.55&   1.55E-05   &1.98   \\
1676.06&   3.88E-06   &2.00   \\\Xhline{1pt}
\end{tabular}
\begin{tabular}{ |c|  c |c|c|}
\Xhline{1pt}
\multicolumn{4}{ |c |}{Test 3 ($\alpha = 1,\beta = 0$)}	\\\cline{1-4}
 R	&time(s) & max error & order\\\hline
1.92E-02   & 0.07&    7.18E-02 & --        \\
 3.19E-03   &0.13&    3.08E-02 & 1.17  \\
8.39E-04    &1.23&    1.10E-02 & 1.40  \\
 1.55E-04  &10.39&    4.14E-03 & 1.45      \\
 2.00E-04 &130.04&    1.23E-03 & 1.75   \\
  8.08E-05  &2052.19&    -- & --    \\\Xhline{1pt}
\end{tabular}
\begin{tabular}{ |c|c|c|c|}
\Xhline{1pt}
\multicolumn{4}{ |c |}{Test 4 ($\alpha = 10,\beta = 0$)}	\\\cline{1-4}
iteration&time(s) 	& max error & order\\\hline
20& 0.85   & 2.81E-01 & --        \\
22& 2.27  & 2.20E-01 & 0.33  \\
21&17.02   & 1.30E-01 & 0.72  \\
20& 173.45  & 5.39E-02 & 1.31      \\
21&2331.48 & 1.39E-02 & 1.95   \\
21&  93834.32  & -- & --    \\\Xhline{1pt}
\end{tabular}
\end{small}
}
\caption{The number of iterations, CPU time, error of numerical solution and order with respect to $\mathcal{N}$ for the transport boundary case.}
\label{tab:fdm_T}
\end{table}

\subsection{RBM results for the parameterized transport boundary problem}

In this section, we present numerical results on the following five problems to demonstrate the applicability and the efficiency of the R2-ROC method for the parameterized transport boundary problem of the Monge-Amp$\grave{\rm e}$re equation.
\paragraph{RB-Test 1}Transporting the following density  to a uniform density on the square $(-0.5,0.5)^2$
\[
f_{X}(\pmb{x},\mu)=1-\frac{0.031}{\mu},
\]
where the auxiliary function is given by $q_\mu(z) = (-\frac{1}{\mu\pi}z^2+\frac{1}{32 \mu \pi ^{3}}+\frac{1}{4 \mu \pi})\cos(8\pi z)+\frac{1}{4 \mu\pi^{2}}z\sin(8\pi z)$ and the exact solution is provided as
$\nabla u = \left(x_{1}+4q'(x_{1},\mu)q(x_{2},\mu),\, x_{2}+4q'(x_{2},\mu)q(x_{1},\mu)\right).$

\paragraph{RB-Test 2}Transporting a uniform density  to a density that blows up at a moving point $(\mu_1, \mu_2)$ on the square $(0,1)^2$:
\[
f_{Y}(\pmb{y},\pmb{\mu})= \frac{\exp(-2\sqrt{(y_{1}-0.5)^2+(y_{2}-0.5)^2})}{\sqrt{(y_{1}-\mu_{1})^2+(y_{2}-\mu_{2})^2}}.
\]

\paragraph{RB-Test 3}Transporting a uniform density to the following density function on the square $(0,1)^2$
\[
\begin{aligned}
f_{Y}(\pmb{y},\mu)&= 1+5\exp(-50|(y_{1}-0.5-\mu)^{2}+(y_{2}-0.5)^{2}-0.09|).
\end{aligned}
\]

\paragraph{RB-Test 4}Transporting a uniform density to the following density function on the square $(0,1)^2$
\[
\begin{aligned}
f_{Y}(\pmb{y},\mu)&= 1+5\exp(-50|(y_{1}-0.5-0.25\cos(2\pi\mu))^{2}+(y_{2}-0.5-0.25\cos(2\pi\mu))^{2}-0.01|).
\end{aligned}
\]

\paragraph{RB-Test 5}Transporting a uniform density on the square $(-0.5,0.5)^2$ to the following density on the disk $y_1^2+y_2^2\le 0.5^2$:
\[
f_{Y}(\pmb{y},\mu)= 1+\exp(-\frac{y_1^2+y_2^2}{2\mu^2})/ (2\pi\mu^2).
\]

\begin{figure}[htbp]
\centering
\subfigure[$\mu = 5$]{
\includegraphics[width=0.21\textwidth]{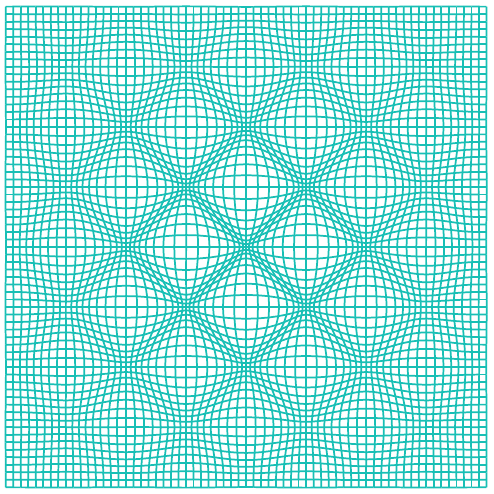}
}
  \subfigure[$\mu = 8$]{
\includegraphics[width=0.21\textwidth]{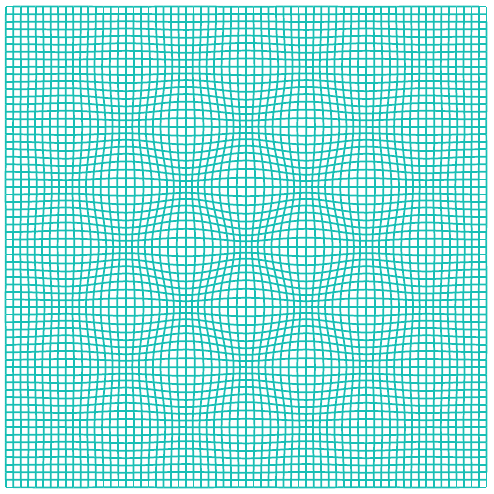}
}
\subfigure[$\mu = 10$]{
\includegraphics[width=0.21\textwidth]{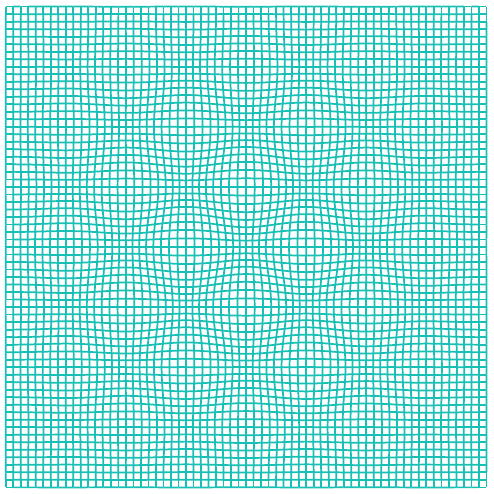}
}
\subfigure[$\mu = 20$]{
\includegraphics[width=0.21\textwidth]{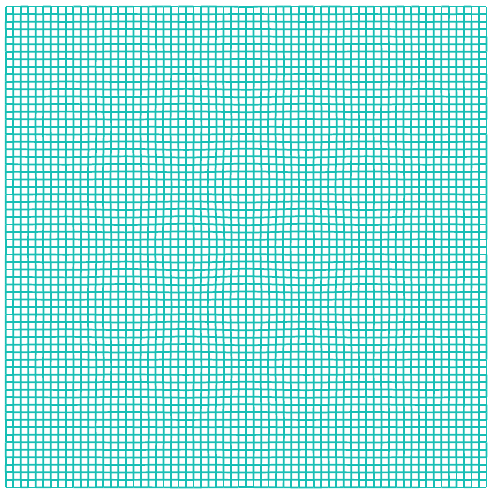}
}\\
\subfigure[$\pmb{\mu} = (0.3,\ 0.6)$]{
\includegraphics[width=0.21\textwidth]{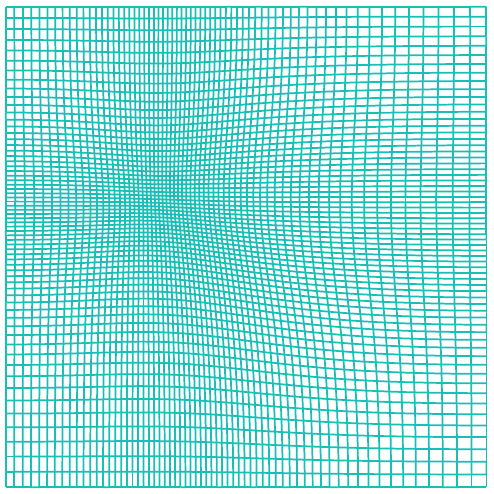}
}
  \subfigure[$\pmb{\mu} = (0.1,\ 0.9)$]{
\includegraphics[width=0.21\textwidth]{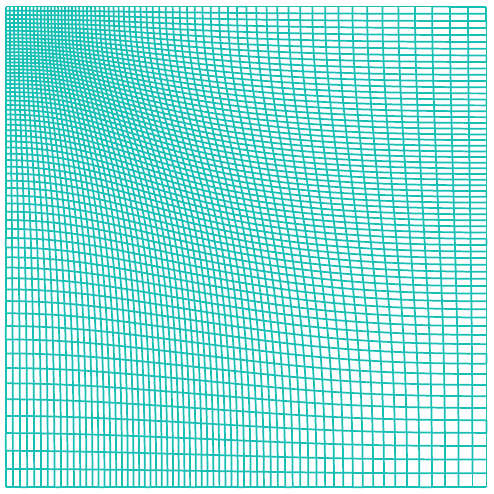}
}
\subfigure[$\pmb{\mu} = (0.5,\ 0.5)$]{
\includegraphics[width=0.21\textwidth]{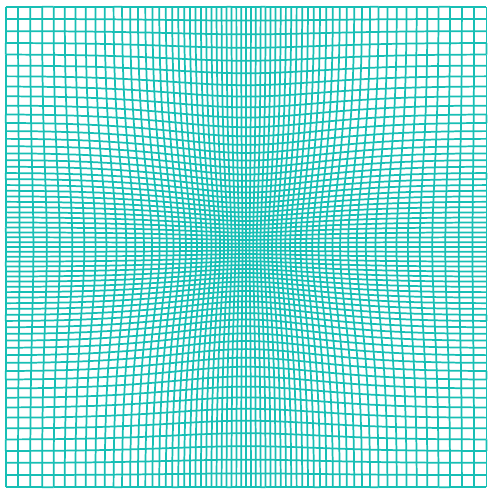}
}
\subfigure[$\pmb{\mu} = (0.8,\ 0.7)$]{
\includegraphics[width=0.21\textwidth]{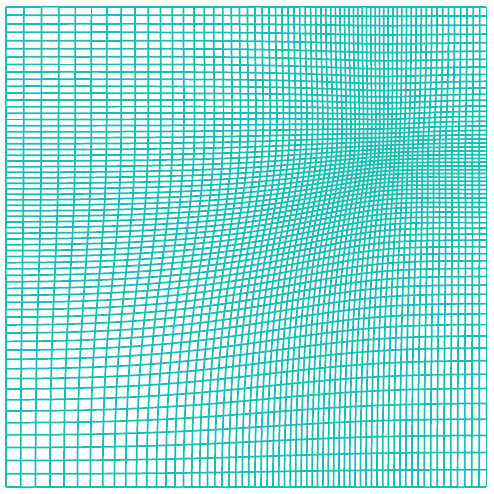}
}\\
\subfigure[$\mu = 0$]{
\includegraphics[width=0.21\textwidth]{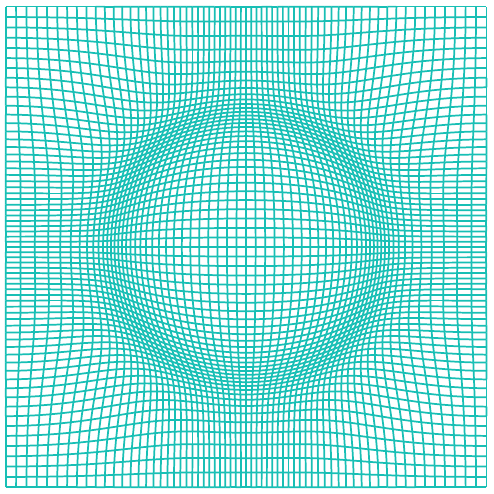}
}
  \subfigure[$\mu = 0.25$]{
\includegraphics[width=0.21\textwidth]{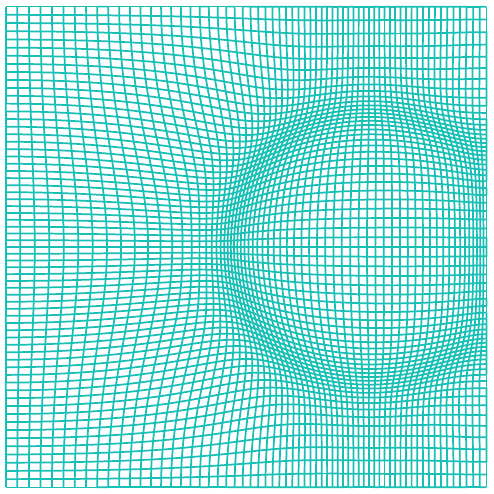}
}
\subfigure[$\mu = 0.5$]{
\includegraphics[width=0.21\textwidth]{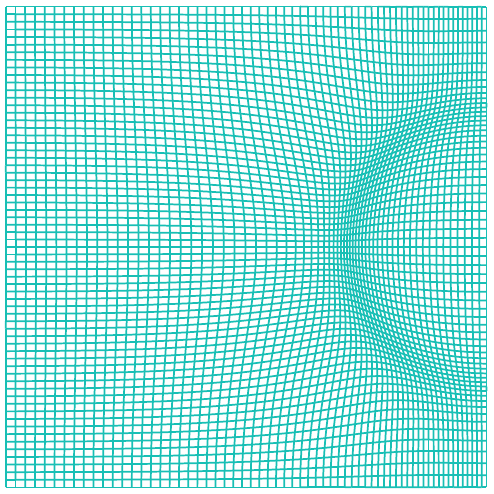}
}
\subfigure[$\mu = 0.75$]{
\includegraphics[width=0.21\textwidth]{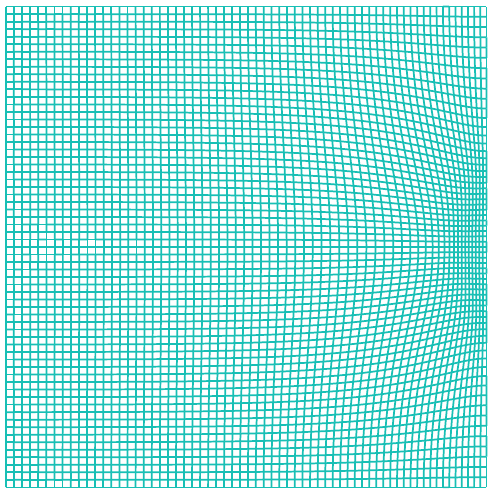}
}\\
\subfigure[$\mu = 0$]{
\includegraphics[width=0.21\textwidth]{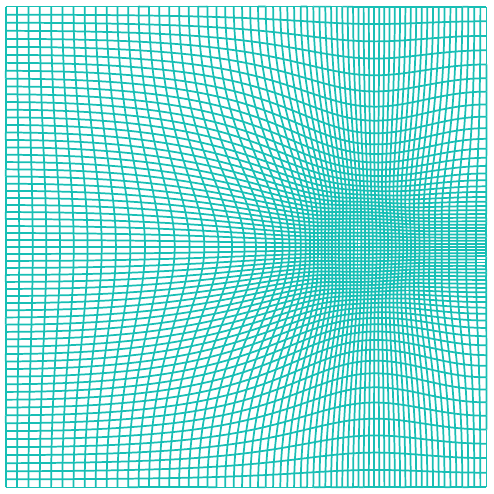}
}
  \subfigure[$\mu = 0.25$]{
\includegraphics[width=0.21\textwidth]{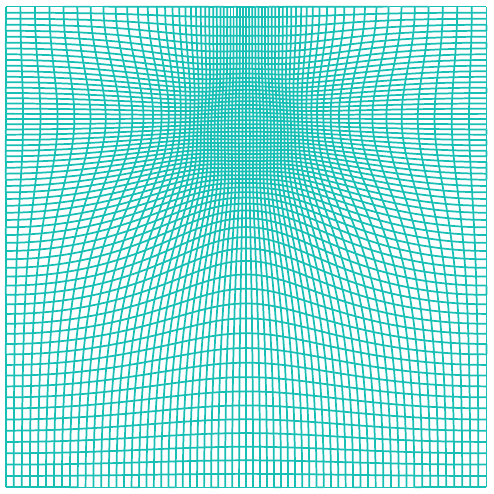}
}
\subfigure[$\mu = 0.5$]{
\includegraphics[width=0.21\textwidth]{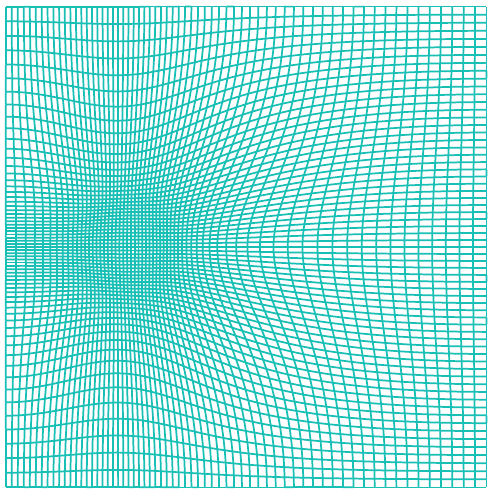}
}
\subfigure[$\mu = 0.75$]{
\includegraphics[width=0.21\textwidth]{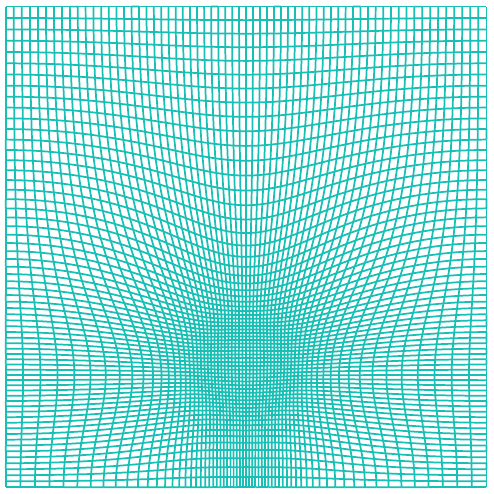}
}\\
\subfigure[$\mu = 0.1$]{
\includegraphics[width=0.21\textwidth]{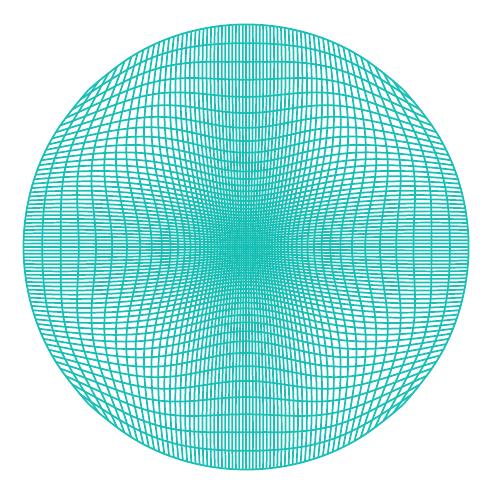}
}
  \subfigure[$\mu = 0.14$]{
\includegraphics[width=0.21\textwidth]{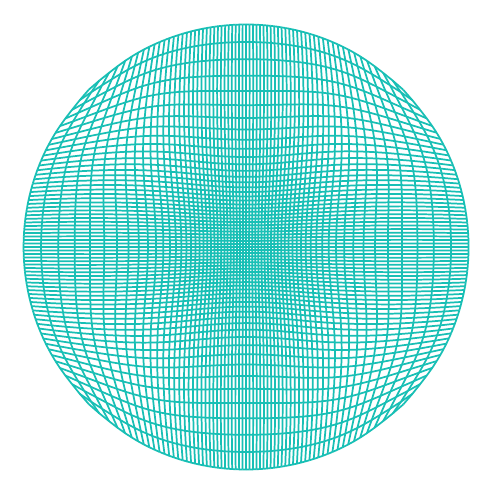}
}
\subfigure[$\mu = 0.19$]{
\includegraphics[width=0.21\textwidth]{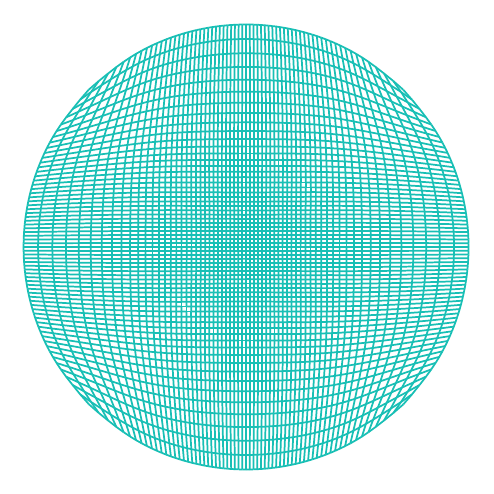}
}
\subfigure[$\mu = 0.24$]{
\includegraphics[width=0.21\textwidth]{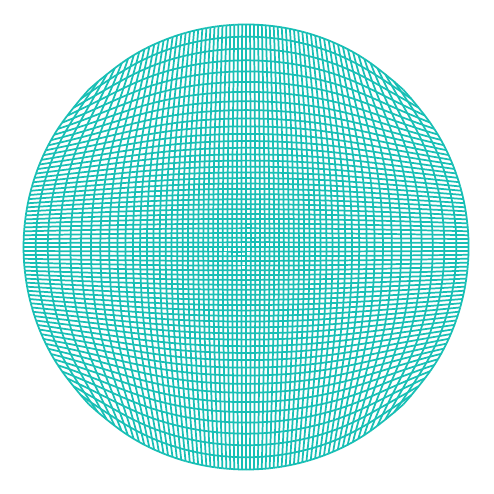}
}\\
\caption{The image of the truth solutions at representative parameter values for test 1 to test 5 of Table \ref{tab:fdm_test_t_roc} (from Top to Bottom).}
\label{fig:para_t_solns}
\end{figure}
\noindent Figure \ref{fig:para_t_solns} shows the truth approximations at representative parameter values that are generated with the narrow-stencil FDM on a mesh of size $\mathcal{N} = 127^2$ and $(\alpha, \beta)$ values given in Table \ref{tab:fdm_test_t_roc}. Parametric variations are clearly visible for each example. In particular, we aim to capture a moving singular point for RB-Test 2 and a circle of denser measure moving to the right for RB-Test 3.

\begin{figure}[htbp]
\centering
\subfigure{
\includegraphics[width=0.31\textwidth]{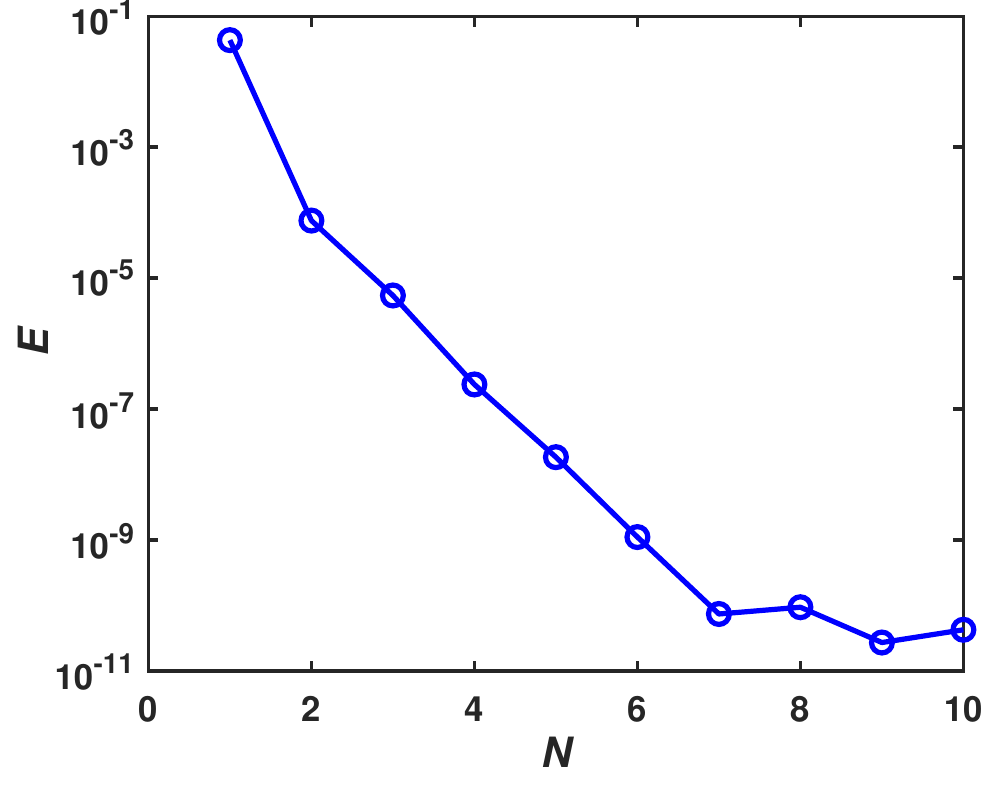}
}
\subfigure{
\includegraphics[width=0.31\textwidth]{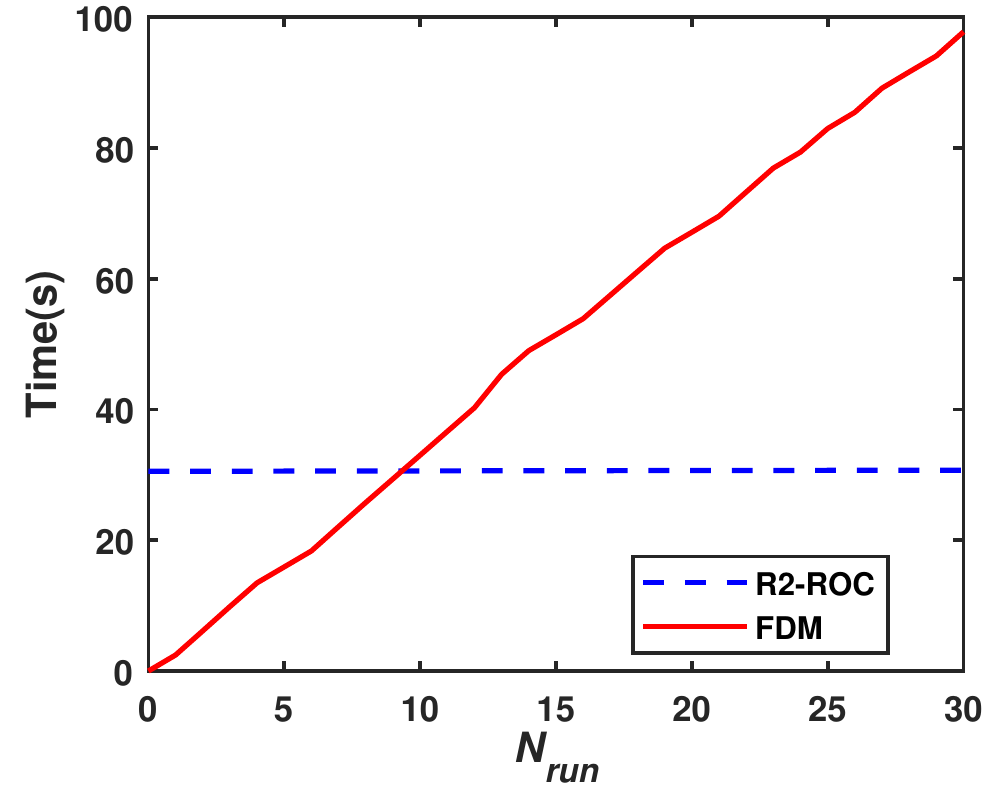}
}
\subfigure{
\includegraphics[width=0.31\textwidth]{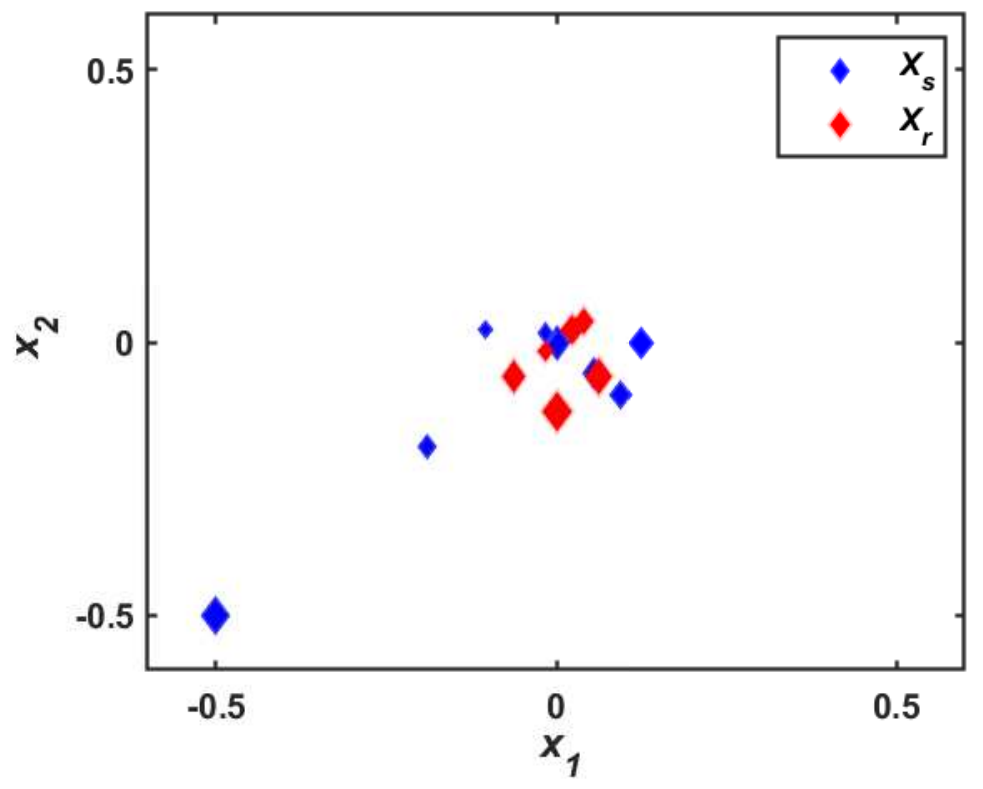}
}
\\
  \subfigure{
\includegraphics[width=0.31\textwidth]{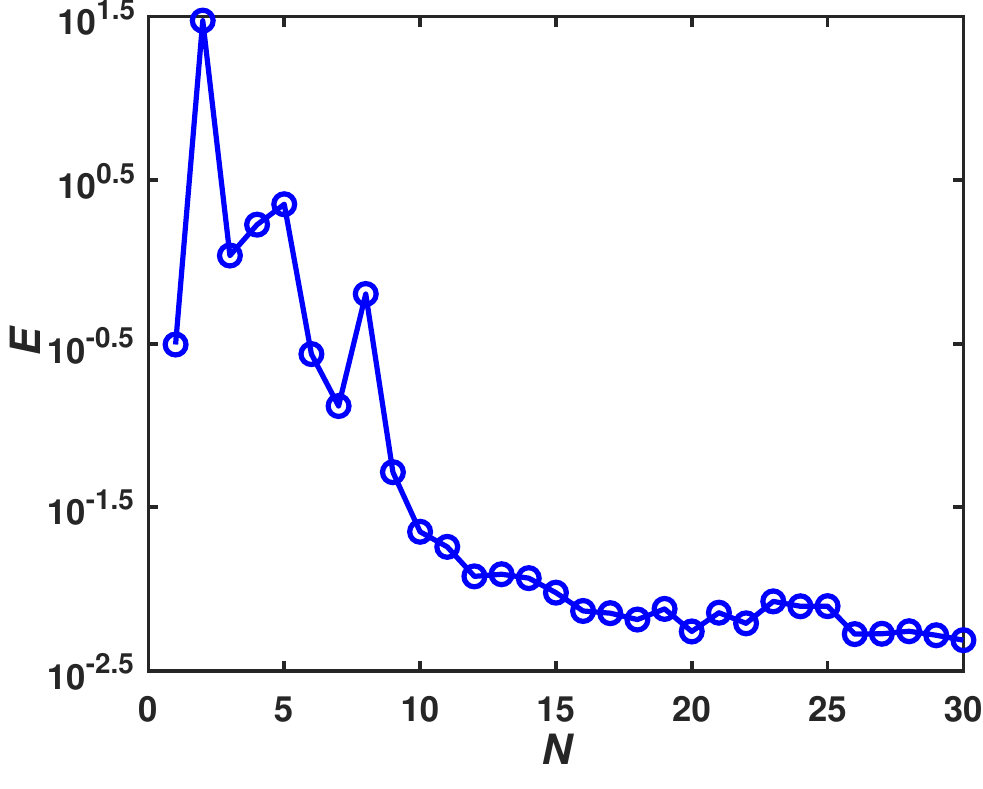}
}
  \subfigure{
\includegraphics[width=0.31\textwidth]{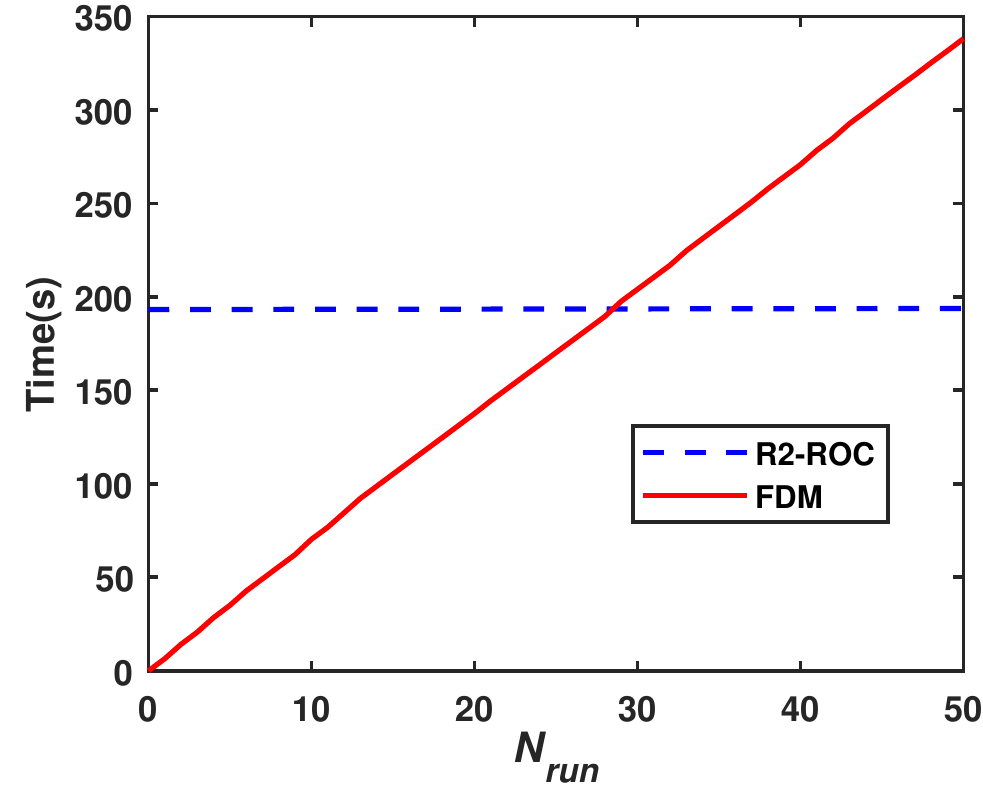}
}
  \subfigure{
\includegraphics[width=0.31\textwidth]{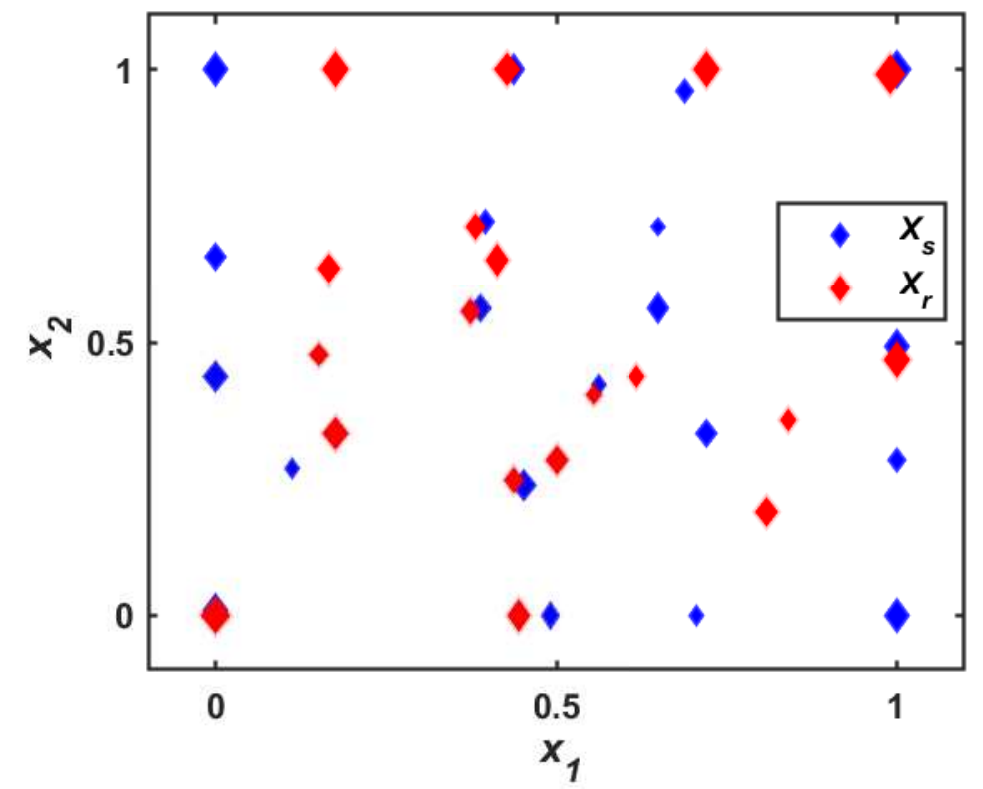}
}
\\
\subfigure{
\includegraphics[width=0.31\textwidth]{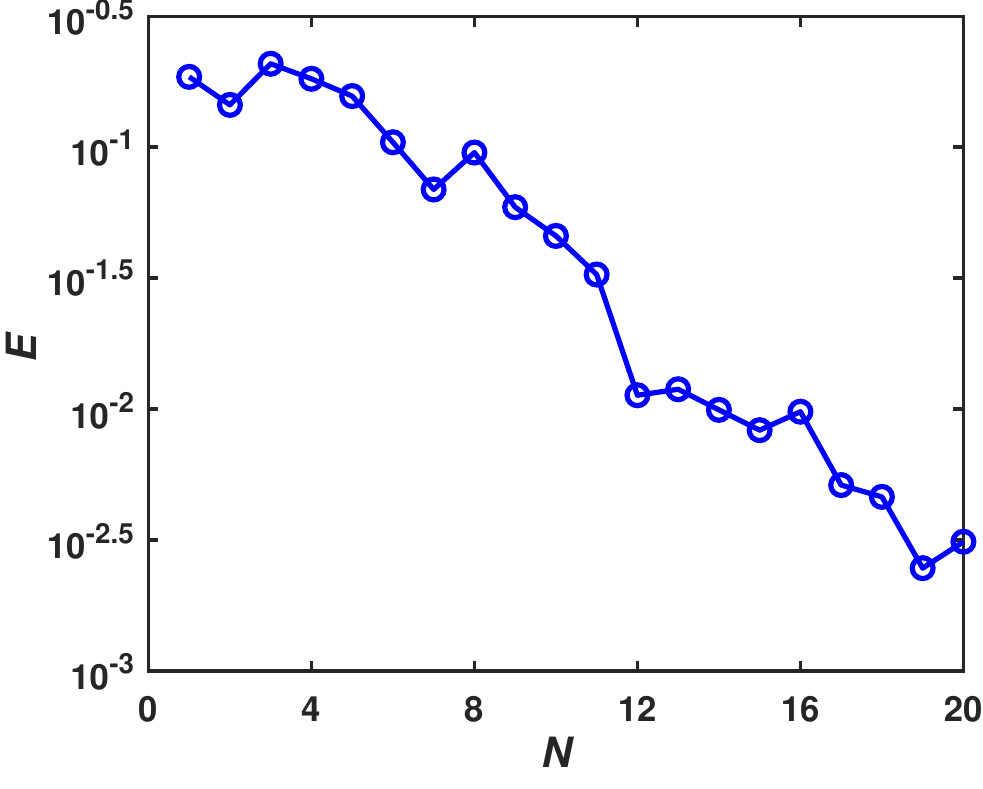}
}
\subfigure{
\includegraphics[width=0.31\textwidth]{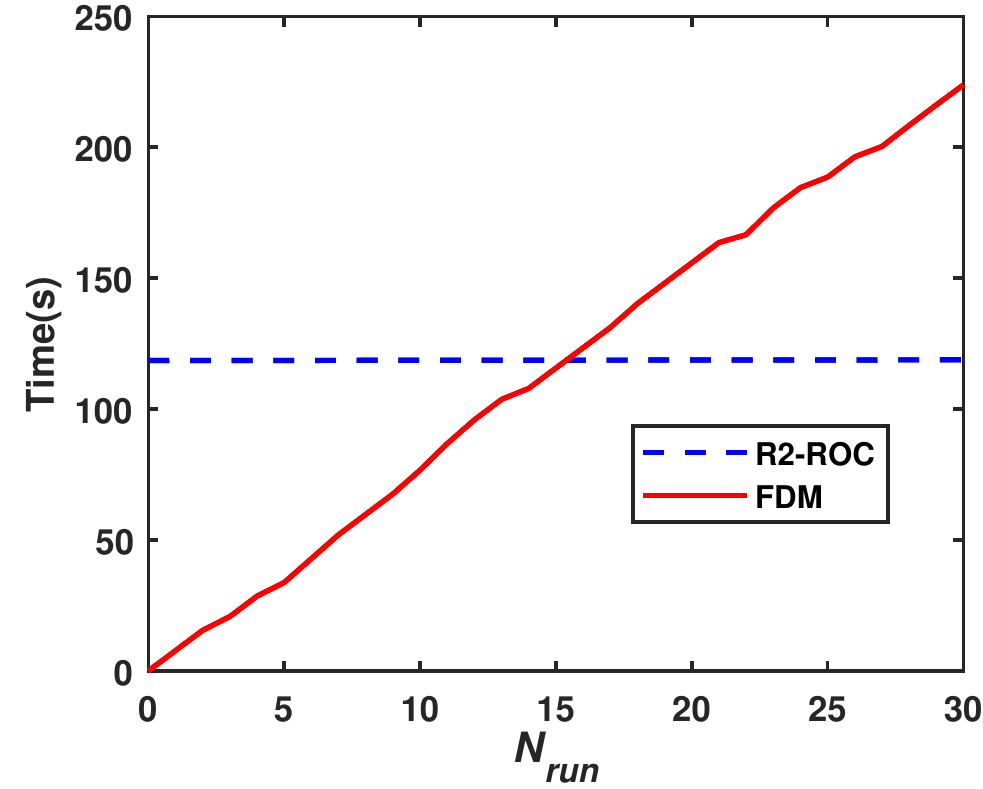}
}
\subfigure{
\includegraphics[width=0.31\textwidth]{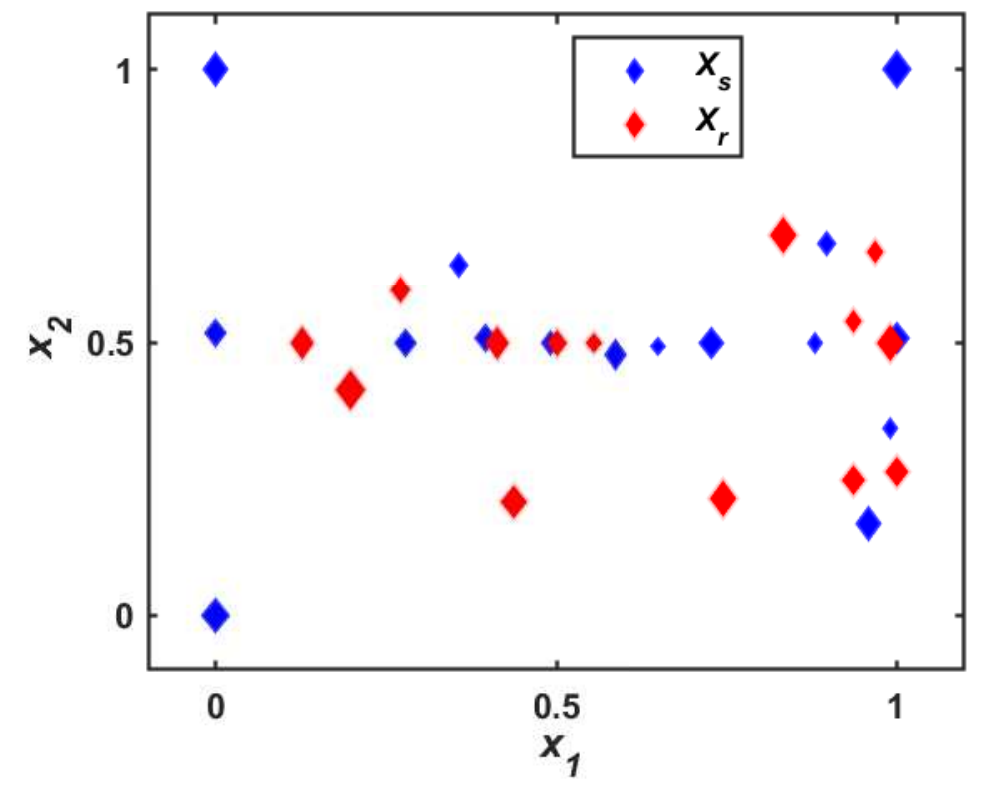}
}
\\
\subfigure{
\includegraphics[width=0.31\textwidth]{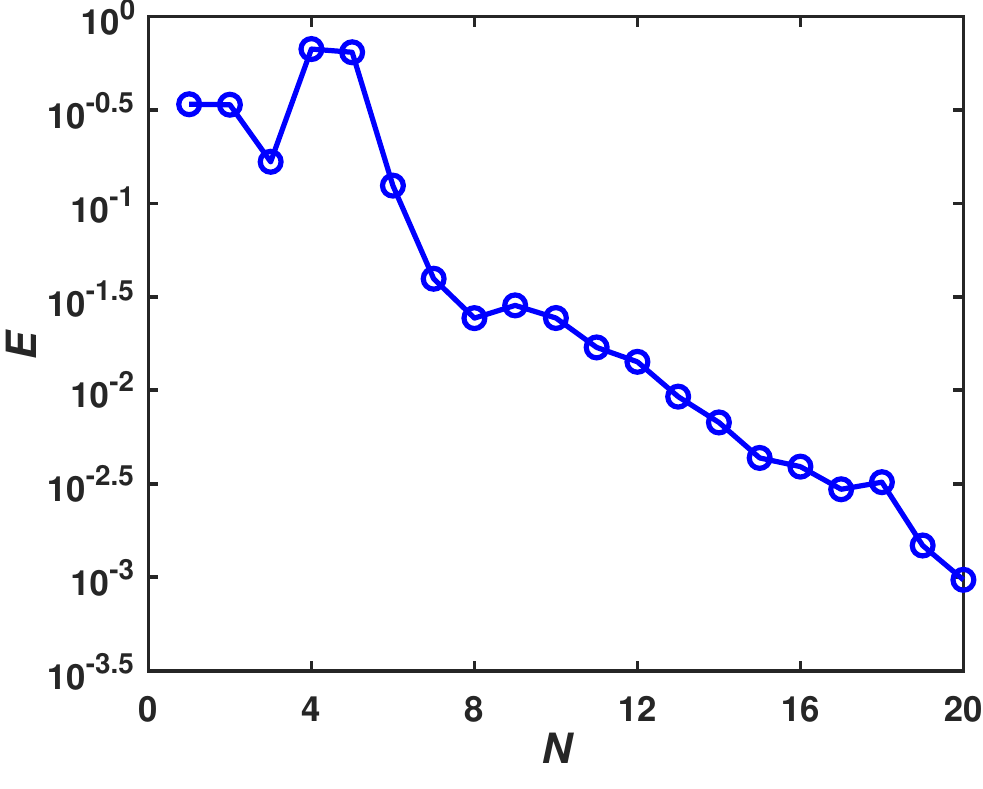}
}
\subfigure{
\includegraphics[width=0.31\textwidth]{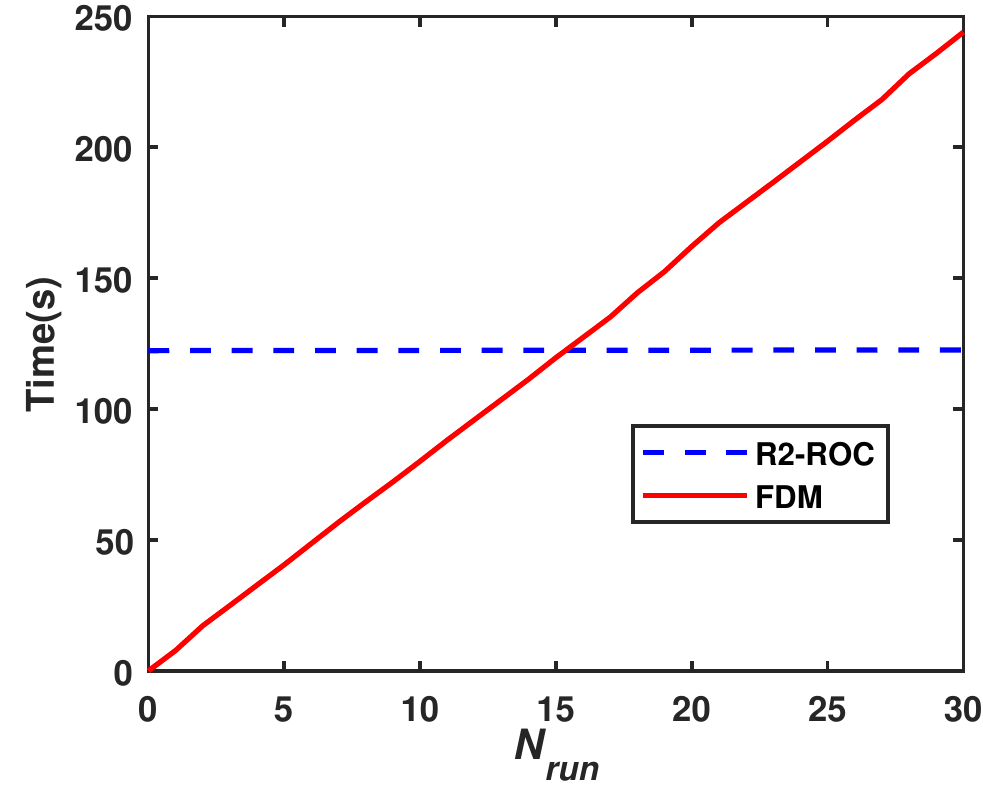}
}
\subfigure{
\includegraphics[width=0.31\textwidth]{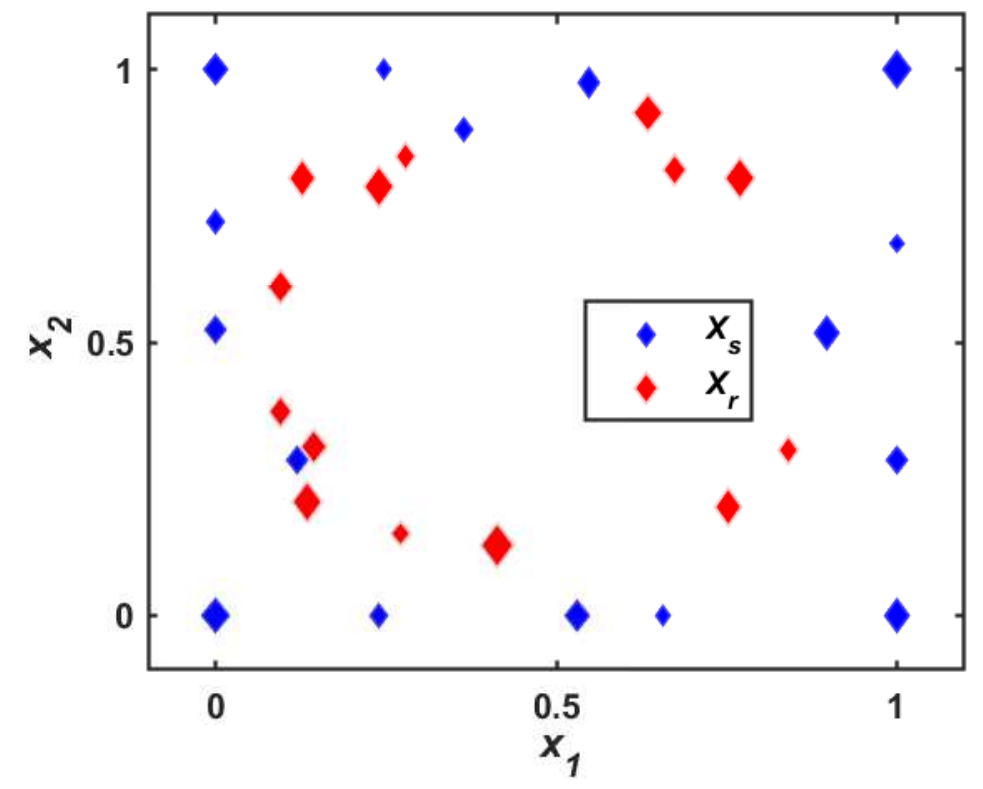}
}
\\
\subfigure{
\includegraphics[width=0.31\textwidth]{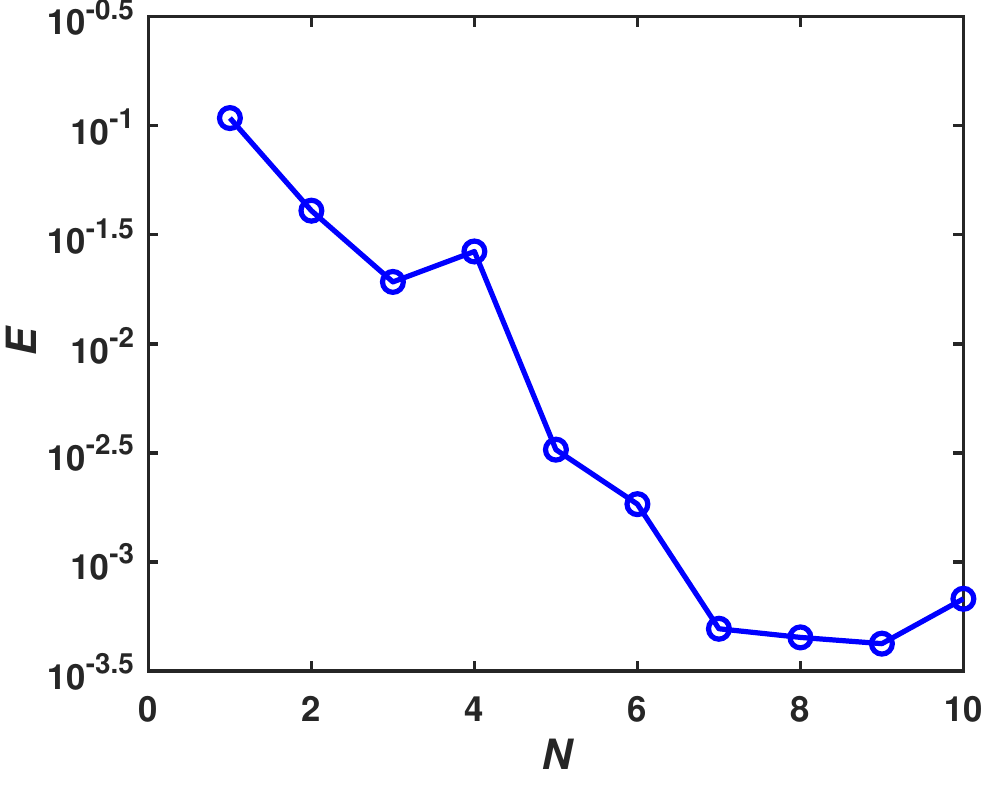}
}
\subfigure{
\includegraphics[width=0.31\textwidth]{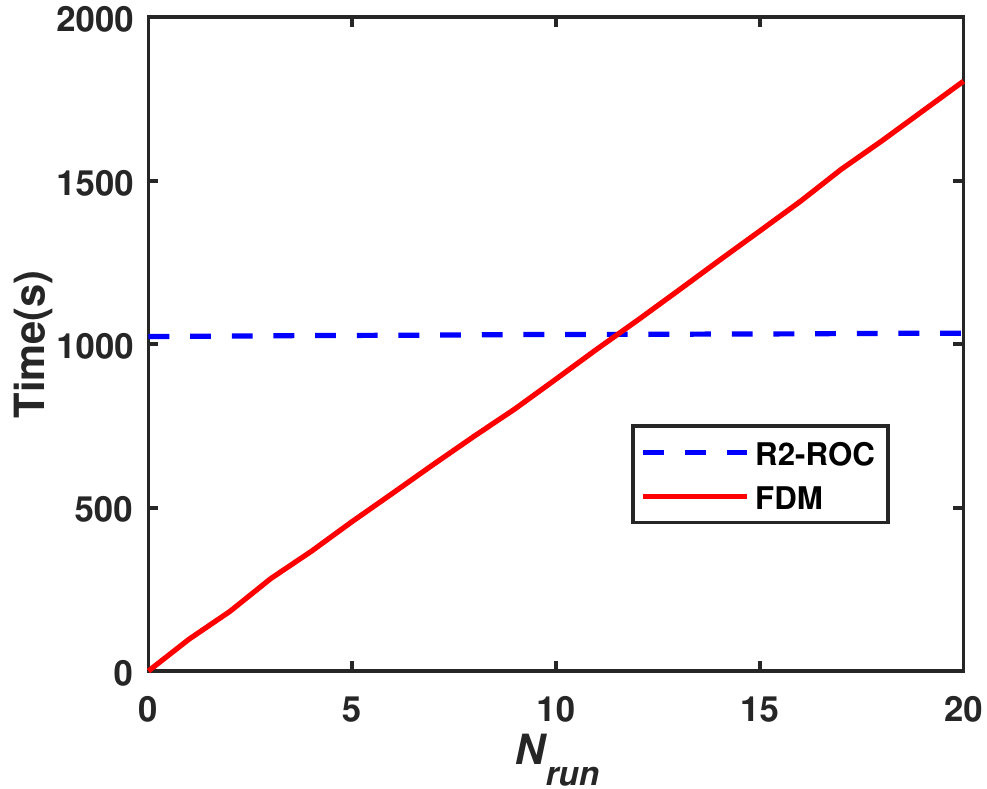}
}
\subfigure{
\includegraphics[width=0.31\textwidth]{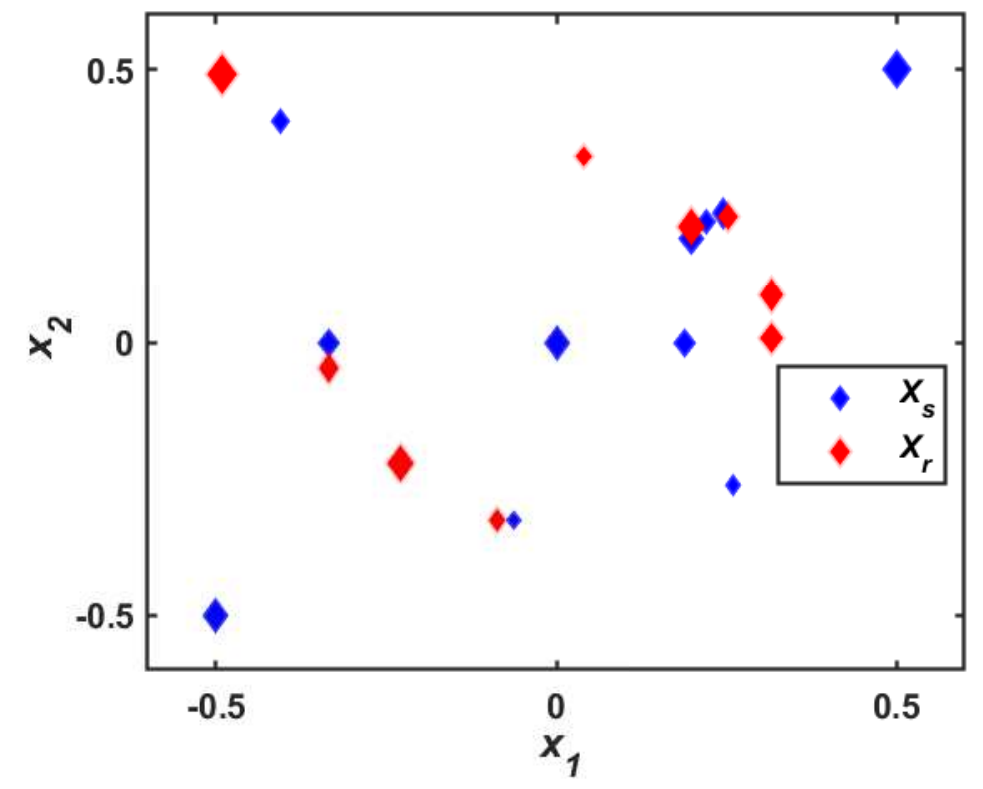}
}
\caption{
R2-ROC results for the parameterized transport boundary problem of the Monge-Amp$\grave{\rm e}$re equation: The histories of convergence (Left), comparison in cumulative run time (Middle), and the collocation points for the R2-ROC method (Right). On the top is for the test 1 and the bottom for the test 5 of Table \ref{tab:fdm_test_t_roc}.
}
\label{fig:fdm_test_t_roc_12}
\end{figure}

\begin{table}[htb]
\centering
\renewcommand{\arraystretch}{1.5}
\begin{tabular}{|c|c|c|c|}
\Xhline{1pt}
RB-Test& $(\alpha, \beta)$ & $\Xi_{\rm train}$ & $\Xi_{\rm test}$ \\
\hline
1& $(0,0)$ & $(5:0.2:20)$ & $(5.1:0.2:19.9)$\\
2& $(200,0)$ & $(0.1:0.04:0.9)^2$ &  $(0.13:0.08:0.89)^2$\\
3& $(50,0)$ & $(0:0.02:1)$ & $(0.01:0.02:9)$\\
4& $(50,0)$ & $(0:0.02:1)$ &  $(0.01:0.02:9)$\\
5&$(10,0)$ & $(0.1:0.01:0.3)$ &  $(0.105:0.01:0.295)$\\
\Xhline{1pt}
\end{tabular}
    \caption{Test problems setup for the parametric transport boundary case.}
    \label{tab:fdm_test_t_roc}
\end{table}
{For these calculations, the error tolerance $\epsilon$ is $10^{-8}$, the maximum number of iterations $K$ is 100.
Using the training and testing sets specified in Table \ref{tab:fdm_test_t_roc} for the R2-ROC method, we generate the reduced basis space and the collocation set, with which we compute the RB solution $\widehat{u}_{N}(\pmb{\mu})$, where $N$ is the number of basis functions that we used. To test the R2-ROC method, we compute the maximum error $E(N)$ between the mappings induced by the RB solution $\widehat{u}_{N}(\pmb{\mu})$ and the truth approximation $u^{\mathcal{N}}(\pmb{\mu})$ for all $\pmb{\mu}\in \Xi_{\rm test}$}. That is, 
\[
E(N) = \max_{\pmb{\mu}\in \Xi_{\rm test}}\left \{||\nabla_{h}u^{\mathcal{N}}(\pmb{\mu})-\nabla_{h}\widehat{u}_{N}(\pmb{\mu})||_{L^{\infty}(X^{\mathcal{N}})}\right\}.
\]
{The left column of Figures \ref{fig:fdm_test_t_roc_12} shows the maximum errors that are plotted against the number of basis functions. Exponential convergence is observed for all cases. The three tests RB-Test 2 to 4 are challenging due to  moving singularities (i.e. the regions of low regularity vary with the parameter $\pmb{\mu}$). The convergence is noticeably slower. The same can be seen for the RB-Test 5 featuring a more challenging transport from a rectangular to a non-rectangular target domain. However, it is clearly still worthwhile to invest in the offline process of the R2-ROC for all cases even when only a modest number of inquiries are needed. To demonstrate the efficiency, we compute the cumulative run time as a function of the number of queries $N_{\rm run}$ for the full and reduced solver with $N$ basis, plotted in the middle column of Figures \ref{fig:fdm_test_t_roc_12}. The offline cost is counted as an overhead for the reduced solver. We then evaluate the break-even number of queries $N_{\rm run}^{\rm e}$ above which it is more costly to run the full simulations for each query (and thus worthwhile to invest the overhead cost training the reduced solver). These quantities and the computation time of the R2-ROC and the full order model are presented in Table \ref{tab:time_t}. As we can see, the break-even numbers of queries $N_{\rm run}^{\rm e}$ are modest across the board thanks to that the computation time of the R2-ROC online solver is hundreds of times smaller than that of the full solver.}

\begin{table}[htb]
\centering
\centering
\begin{small}
\begin{tabular}{ | c | c |c|    c | c|c|}
\Xhline{1pt}
\multicolumn{1}{  |c |}{\multirow{2}{*}{Test}}
&\multicolumn{1}{  c |}{\multirow{2}{*}{$N$}}
&\multicolumn{1}{  c |}{\multirow{2}{*}{$N_{\rm run}^e$}}
&\multicolumn{2}{ c |}{R2-ROC}\\ \cline{4-5}
\multicolumn{1}{ | c |}{}
&&& Offline &Online	\\\hline
1&   7& 10&29.45&   0.0091       \\
2&   20& 30&193.24&   0.011       \\
3&  15& 17&118.53&   0.0093        \\
4&  15& 17&122.32&   0.0088          \\
5&  10& 12&1023.16&   0.5303          \\\Xhline{1pt}
\end{tabular}
\begin{tabular}{ |c|  }
\Xhline{1pt}
\multicolumn{1}{ |c |}{FDM}	\\\cline{1-1}
- 	\\\hline
 3.18         \\
 6.76    \\
7.46    \\
 8.21      \\
  90.24      \\\Xhline{1pt}
\end{tabular}
\end{small}
\caption{Offline and Online computational times for different tests.}
\label{tab:time_t}
\end{table}

\subsection{Dirichlet boundary value problem of the Monge-Amp$\grave{\textbf{e}}$re equation}
For completeness, we test our methods on the Dirichlet boundary value problem of Monge-Amp$\grave{\rm e}$re equation \eqref{eq_MA_D} and report the results in this section. In this case, the narrow-stencil FDM could be applied directly without the need of an iterative procedure for boundary enforcement.
We consider three tests of decreasing regularities. Listed in Table \ref{tab:fdm_test_d} are their density function $f(\pmb{x})$, exact solution $u(\pmb{x})$ (which induces the boundary condition $g(\pmb{x})$), and computational domain $X$. The relative $L^\infty$ error between exact solution $u(\pmb{x})$ and its approximation $u^{\mathcal{N}}(X^{\mathcal{N}})$, and its convergence orders together with the penalization parameters $\alpha$ and $\beta$ are presented in Table \ref{tab:fdm_D}. This verifies that the method converges with the expected order of $2$ when the first derivative is continuous.
\begin{table}[h]
\centering
\renewcommand{\arraystretch}{1.5}
\begin{tabular}{|c|c|c|}
\hline
Test& $(f(\pmb{x}),u(\pmb{x}))$&$X$\\
\hline
\multirow{2}{*}{$C^\infty$} & $f(\pmb{x})=(1+x_{1}^{2}+x_{2}^{2})\exp(x_{1}^{2}+x_{2}^{2})$     & \multirow{2}{*}{$X = (0,1)\times(0,1)$} \\
&$u(\pmb{x})=\exp\left(\frac{x_{1}^{2}+x_{3}^{2}}{2}\right)$     & \\
\hline

\multirow{2}{*}{$C^1$} &$f(\pmb{x})=\left(1-\frac{0.2}{|(x_{1}-0.5,x_{2}-0.5)|}\right)^{+}$     & \multirow{2}{*}{$X = (0,1)\times(0,1)$} \\
&$u(\pmb{x})=\frac{1}{2}((|(x_{1}-0.5,x_{2}-0.5)|-0.2)^{+})^{2}.$   &  \\

\hline
\multirow{2}{*}{$C^0$} & $f(\pmb{x})=0$     & \multirow{2}{*}{$X = (-1,1)\times(-1,1)$} \\
&$u(\pmb{x})=|x_{1}|$    & \\
\hline
\end{tabular}
    \caption{Setup of the Dirichlet test problems.}
    \label{tab:fdm_test_d}
\end{table}
\begin{table}[h]
\centering
\bigskip
\centering
\begin{small}
\begin{tabular}{ | c |   c| c | c|}
\Xhline{1pt}
\multicolumn{1}{ | c |}{\multirow{2}{*}{$\mathcal{N}$}}
&\multicolumn{3}{ c |}{Test $C^\infty$ ($\alpha = 1,\beta = 0$)}\\ \cline{2-4}
\multicolumn{1}{ | c |}{}
&time(s)&error& order	\\\hline
$15^2$&    0.03& 8.98E-03   & --     \\
$31^2$&  0.07 & 2.38E-03   & 1.83 \\
$65^2$&  0.63  &5.97E-04   & 1.87   \\
$127^2$&   4.48  &1.65E-04   & 1.92   \\
$255^2$&  45.74  &4.28E-04   & 1.94 \\
$511^2$&  688.69  &1.09E-05   & 1.96\\\Xhline{1pt}
\end{tabular}
\begin{tabular}{ |c|c|c|}
\Xhline{1pt}
\multicolumn{3}{ |c |}{Test $C^1$ ($\alpha = 10,\beta = 0$)}\\ \cline{1-3}
\multicolumn{1}{ | c |}{}
time(s)&error& order	\\\hline
   0.02& 8.54E-02  &--  \\
 0.05 & 2.07E-02 &1.95   \\
  0.56& 5.38E-03 &1.82    \\
  5.52& 3.05E-03 &0.85 \\
   68.23& 1.70E-03 & 0.84 \\
  1029.31& 8.14E-04 & 1.06 \\\Xhline{1pt}
\end{tabular}
\begin{tabular}{ |  c|c|c|}
\Xhline{1pt}
\multicolumn{3}{ |c |}{Test $C^0$ ($\alpha = 200,\beta = 0$)}\\ \cline{1-3}
\multicolumn{1}{ | c |}{}
time(s)&error& order	\\\hline
  0.02& 1.14E-00  &--  \\
 0.07 & 1.00E-00 &0.19   \\
  0.63& 4.92E-01 &0.96    \\
 4.69& 2.03E-01 &1.32 \\
   80.14& 1.04E-01 & 0.96 \\
   --&--&--\\\Xhline{1pt}
\end{tabular}
\end{small}
\caption{\label{table.label} Computation time, maximum error and rates of convergence for the FDM solutions of the Dirichlet case.}
\label{tab:fdm_D}
\end{table}

\begin{figure}[htbp]
\centering
\subfigure[$\mu = 0.1$]{
\includegraphics[width=0.23\textwidth]{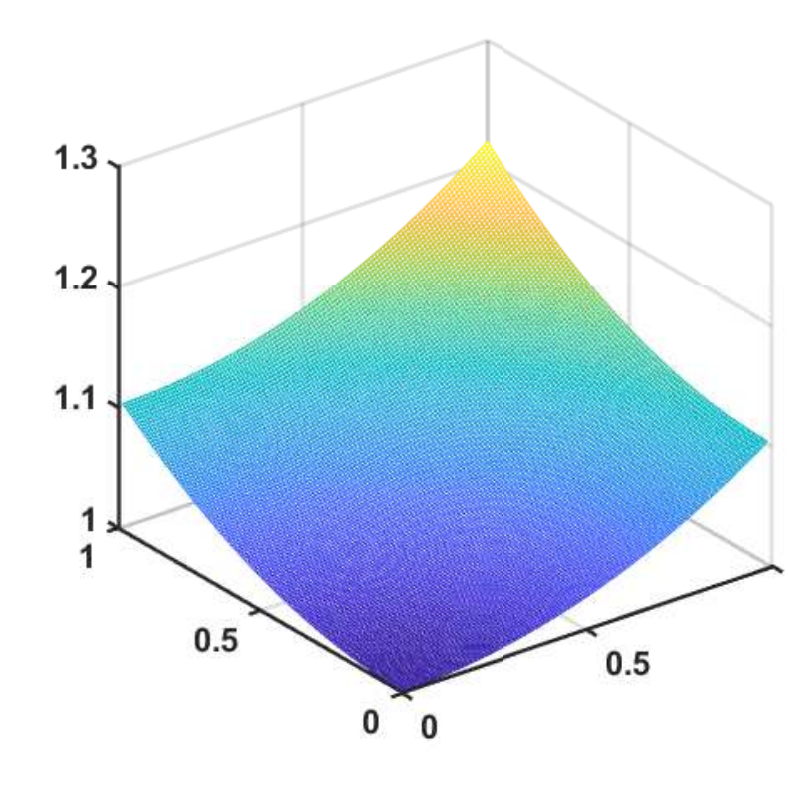}
}
  \subfigure[$\mu = 1$]{
\includegraphics[width=0.23\textwidth]{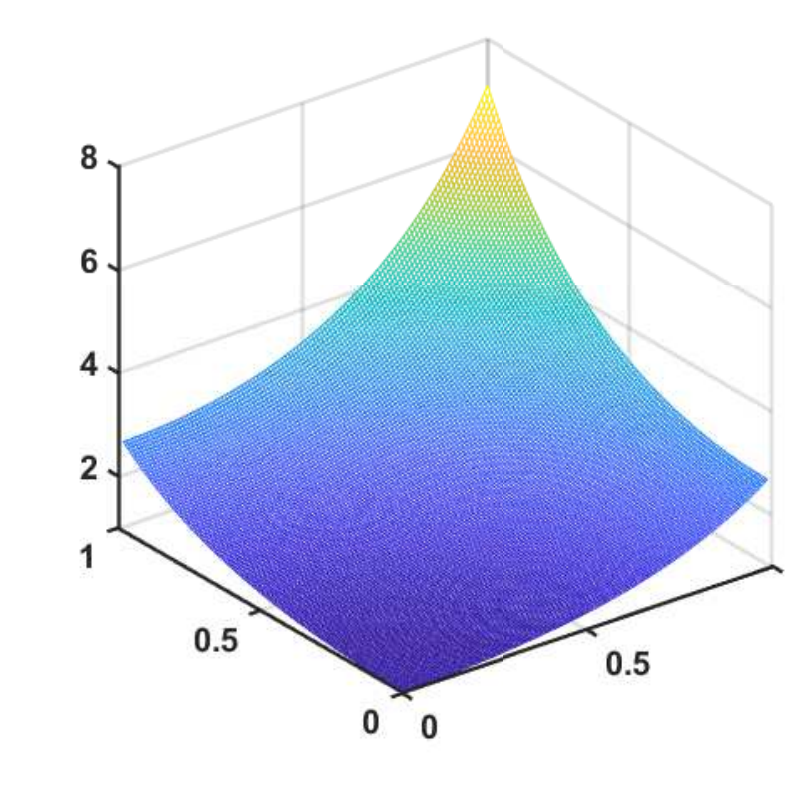}
}
\subfigure[$\mu = 0.1$]{
\includegraphics[width=0.23\textwidth]{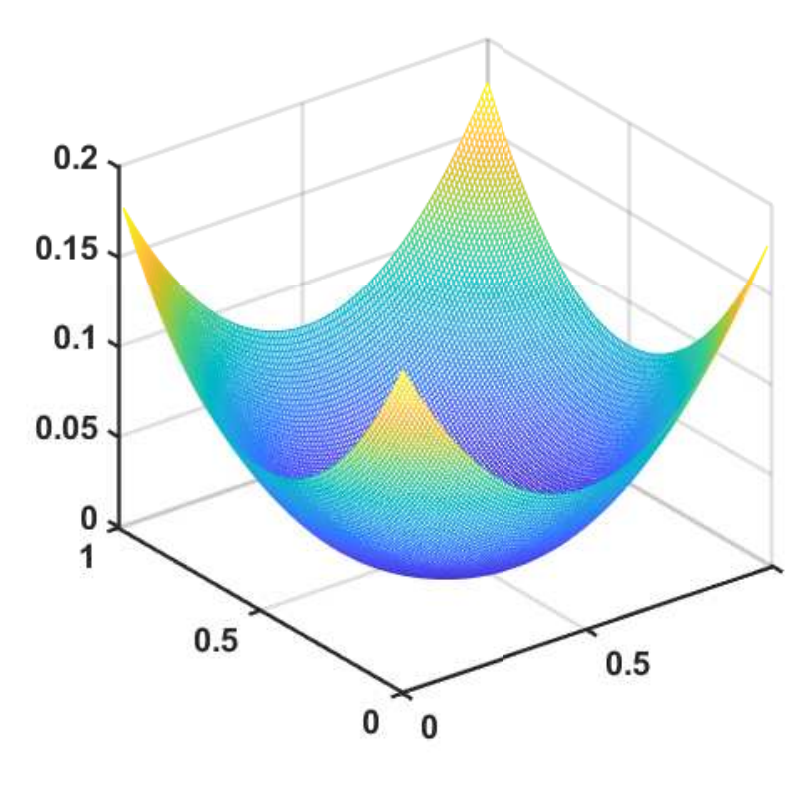}
}
\subfigure[$\mu = 0.5$]{
\includegraphics[width=0.23\textwidth]{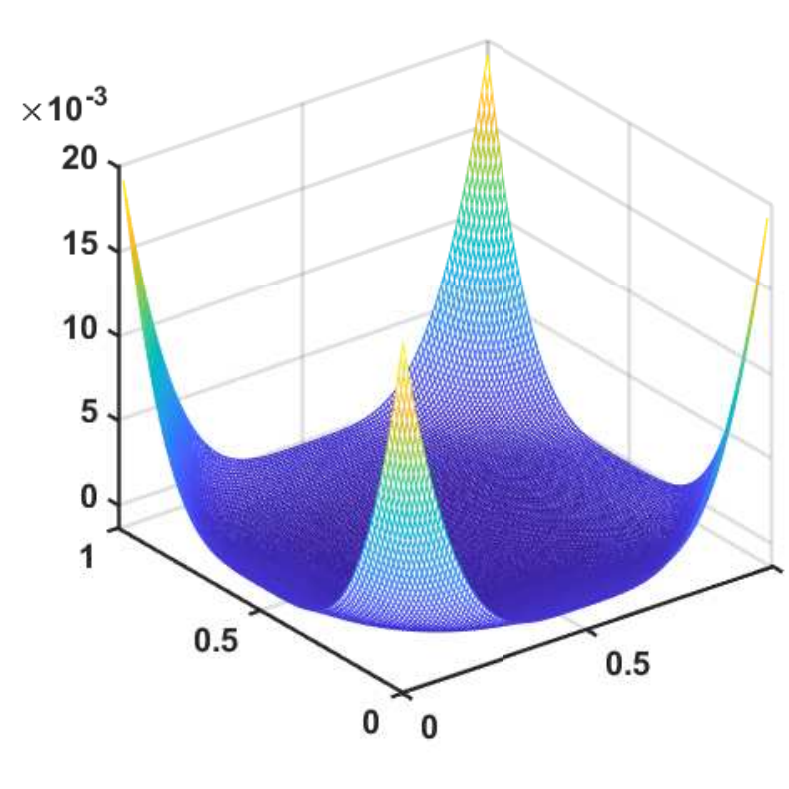}
}
\caption{Truth solutions at representative parameter values for $C^{\infty}(\mu)$ test ((a) and (b)) and $C^1(\mu)$ test ((c) and (d)).}
\label{fig:truth_para_d}
\end{figure}

To test our reduced order solver, we consider two parameterized Monge-Amp$\grave{\rm e}$re equations listed in Table \ref{tab:fdm_test_d_roc}, corresponding to two of the non-parametric cases in Table \ref{tab:fdm_test_d}. Truth approximations are generated with $\alpha = 1$ and $\beta = 0$ on a uniform cartesian mesh of size $\mathcal{N} =127^2$. Figure \ref{fig:truth_para_d} shows two representative solutions for each example. For the first example, solutions for $\mu=0.1$ and $\mu=1$ differ the most around the $(1,1)$-corner of the domain while, for the second example is more challenging with the parameter $\mu$ dictating the location of the regularity change in the solution.

\begin{figure}[htbp]
\centering
\subfigure{
\includegraphics[width=0.31\textwidth]{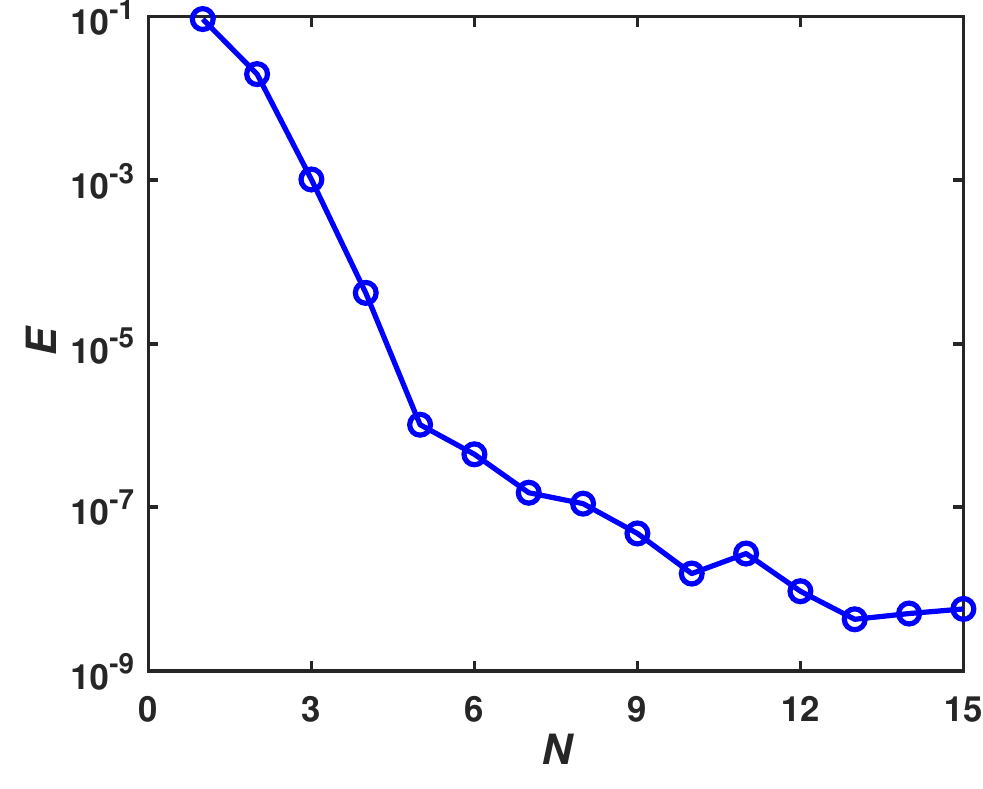}
}
\subfigure{
\includegraphics[width=0.31\textwidth]{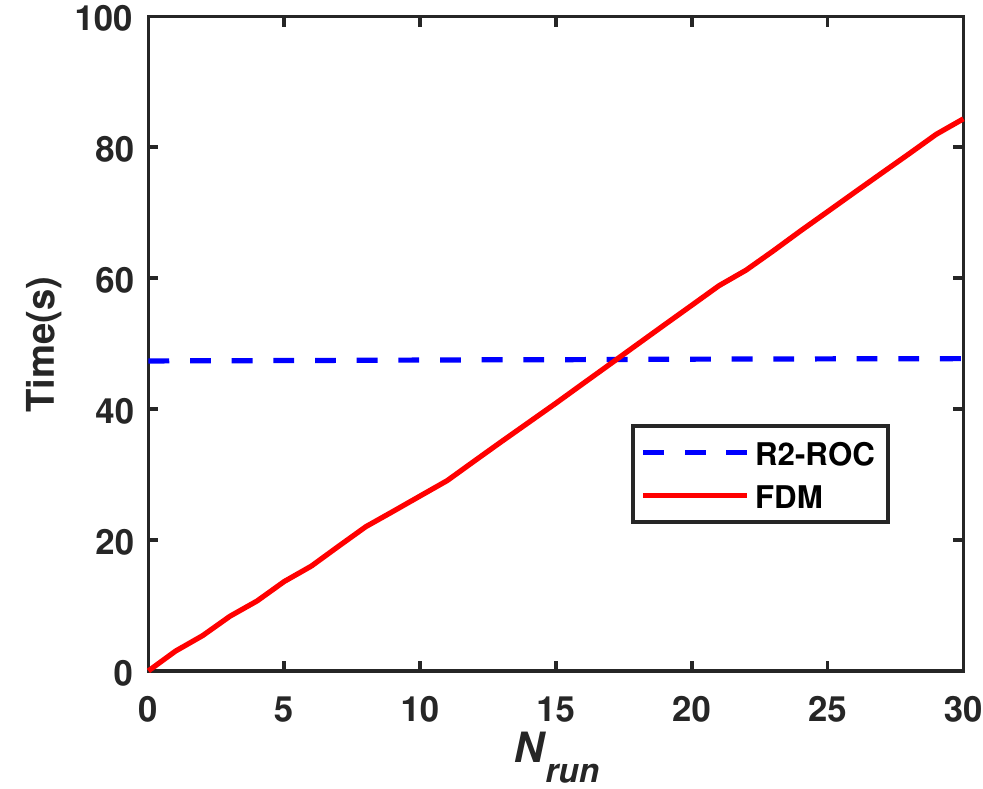}
}
\subfigure{
\includegraphics[width=0.31\textwidth]{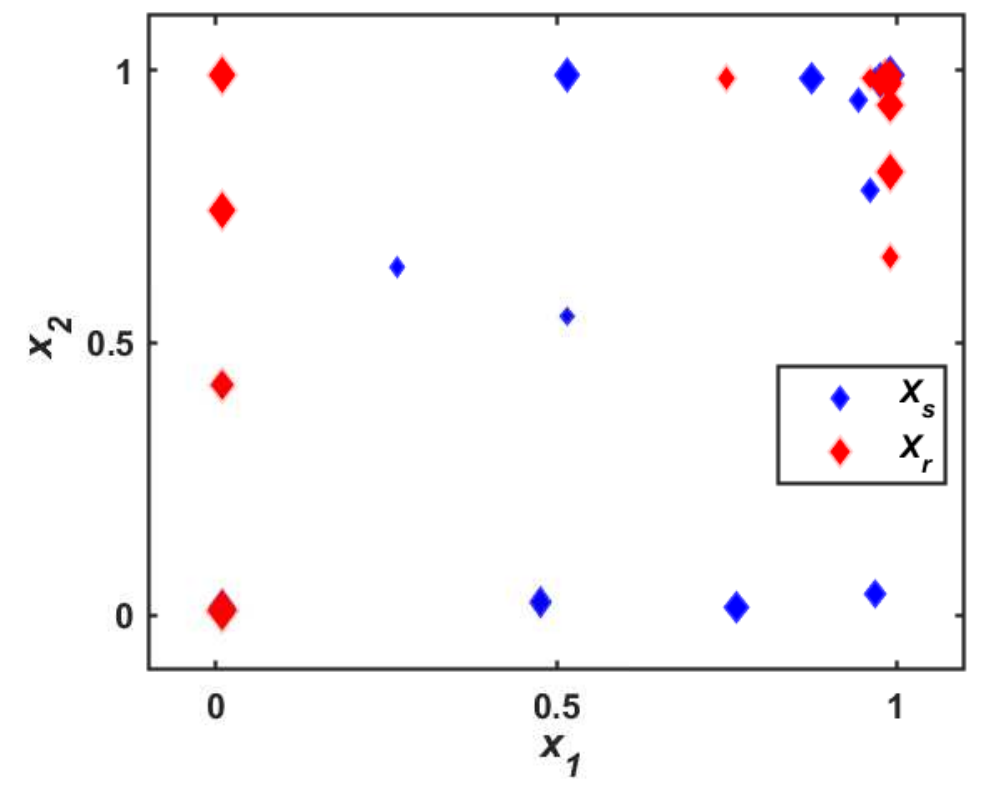}
}
\\
  \subfigure{
\includegraphics[width=0.31\textwidth]{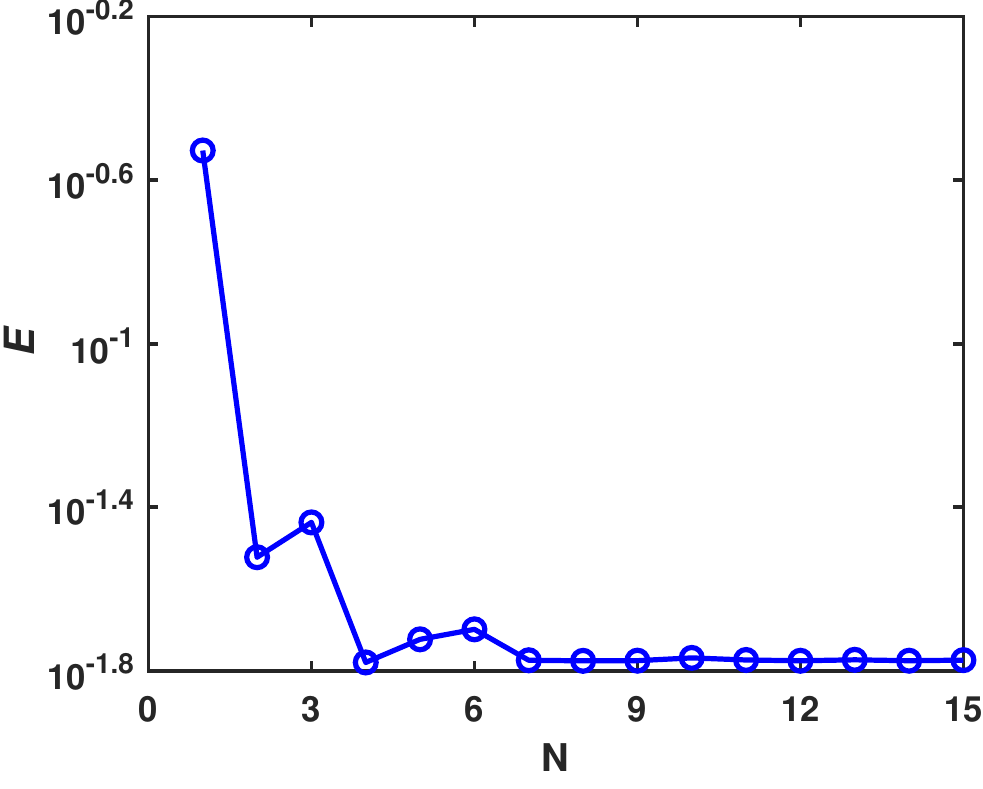}
}
  \subfigure{
\includegraphics[width=0.31\textwidth]{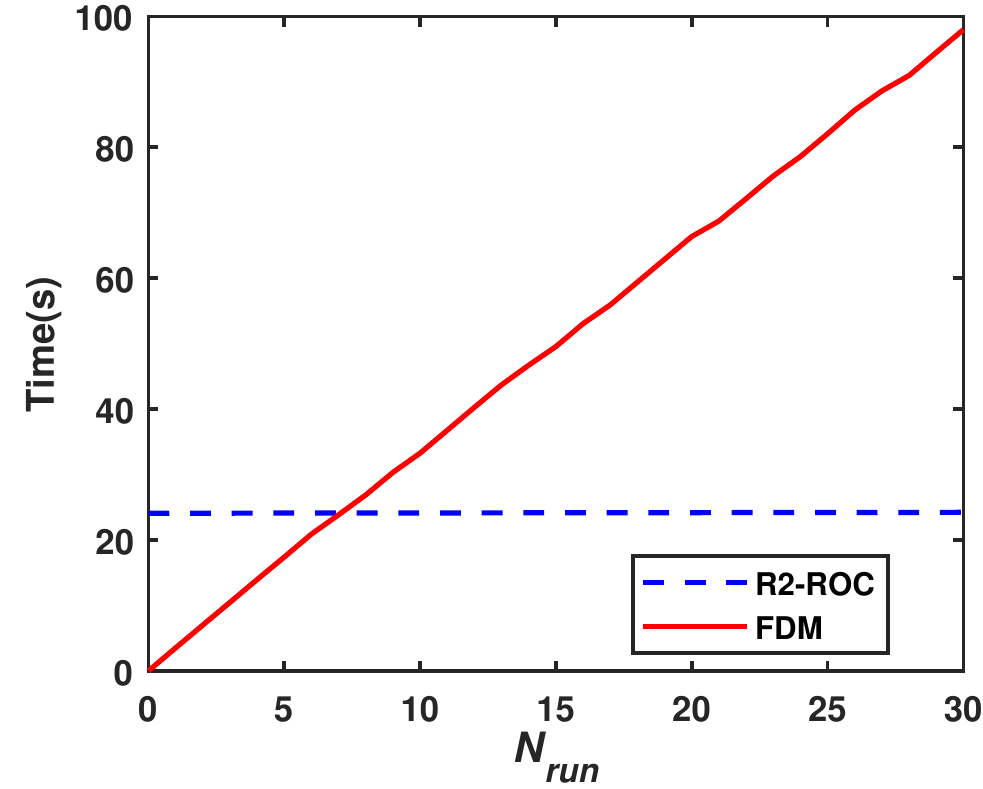}
}
  \subfigure{
\includegraphics[width=0.31\textwidth]{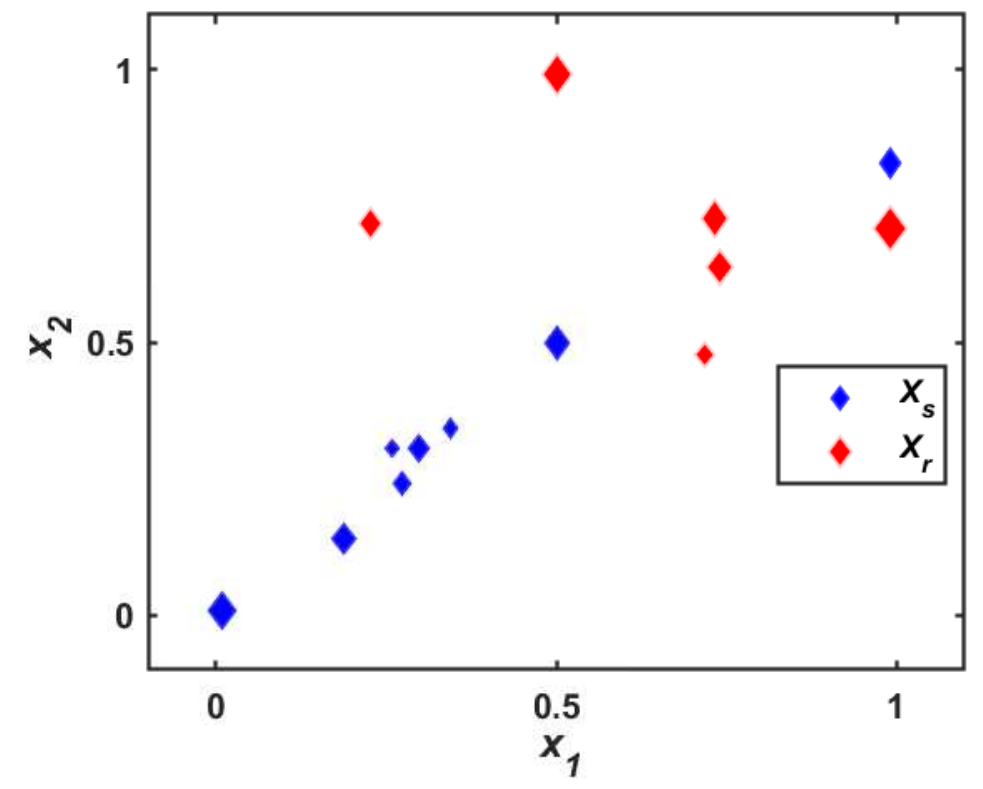}
}
\caption{
R2-ROC results for the parameterized Dirichlet boundary problem of the Monge-Amp$\grave{\rm e}$re equation: The histories of convergence (Left), comparison in cumulative run time (Middle), and the collocation points for the R2-ROC method (Right). On the top is for the $C^{\infty}(\mu)$ test and the bottom for the $C^1(\mu)$ test of Table \ref{tab:fdm_test_d_roc}.
}
\label{fig:r2roc_para_d}
\end{figure}

The R2-ROC results are presented in Figure \ref{fig:r2roc_para_d}. We can see, from the left column, that the errors steadily decrease for the first test while it plateaus for the more challenging second test. The middle column displays the comparison in cumulative computation time for R2-ROC (with $15$ bases for the first example and $7$ for the second) and finite difference method as we increase the number of simulations $N_{\rm run}$. We see that R2-ROC starts to save time when the number of simulations is above $18$ or $8$ for the two tests respectively, and that the savings become more dramatic as  $N_{\rm run}$ increases due to the negligible marginal expense of R2-ROC.
\begin{table}[ht]
\centering
\renewcommand{\arraystretch}{1.5}
\begin{tabular}{|c|c|c|c|}
\hline
RB-Test& $(f(\pmb{x}), u(\pmb{x}))$ & $\Xi_{\rm train}$ & $\Xi_{\rm test}$ \\
\hline
\multirow{2}{*}{$C^\infty (\mu)$} & $f(\pmb{x})= 4\mu^{2}(1+2\mu(x_{1}^{2}+x_{2}^{2}))\exp(2\mu(x_{1}^{2}+x_{2}^{2}))$     & \multirow{2}{*}{$(0.1:0.02:1)$} & \multirow{2}{*}{$(0.11:0.02:0.99)$}\\
&$u(\pmb{x}) = \exp(\mu(x_{1}^{2}+x_{2}^{2}))$   &&   \\
\hline
\multirow{2}{*}{$C^1(\mu)$} &$f(\pmb{x})=\left(1-\frac{\mu}{|(x_{1}-0.5,x_{2}-0.5)|}\right)^{+}$     & \multirow{2}{*}{$(0.1:0.01:.5)$} &  \multirow{2}{*}{$(0.105:0.01:0.495)$}\\
&$u(\pmb{x})=\frac{1}{2}((|(x_{1}-0.5,x_{2}-0.5)|-\mu)^{+})^{2}$  &&   \\
\hline
\end{tabular}
    \caption{Setup of the test problems for the parametric Dirichlet case.}
    \label{tab:fdm_test_d_roc}
\end{table}

\section{Conclusion}
\label{Section:conclusion}
In this article, we develop a fast algorithm for the nonlinear parameterized  Monge-Amp$\grave{\rm e}$re equation with transport boundaries which models the optimal transport problem with $L^{2}$ cost function. It features a novel extension of the narrow-stencil finite difference scheme \cite{feng2019narrow} to our setting incorporating the projection-iteration method \cite{froese2012numerical} to deal with the transport boundary. The resulting solver is shown to be effective and accurate even when facing low-regularity.
Building on this truth approximation solver, we adapt the R2-ROC algorithm \cite{chen2021eim,chen2021l1} to
the parameterized  Monge-Amp$\grave{\rm e}$re equation with transport boundaries. 
Several challenging tests with different parameter delineations are provided to demonstrate the method's capability in efficiently producing an accurate and reliable mapping induced by the RB solution.

\bibliographystyle{abbrv}
\bibliography{ref}

\end{document}